%% file: article.tex
\documentclass[a4paper,12pt]{article}

\pdfoutput=1
\usepackage{amsmath,amssymb,wasysym}

\usepackage[T1]{fontenc}
\usepackage[utf8x]{inputenc}
\usepackage[french]{babel}
\usepackage{enumerate}
\usepackage{url}
\usepackage{pdfpages}
\usepackage{graphicx}
\usepackage{color}

\title{Vénus selon Ibn al-\v{S}\=a\d{t}ir\footnote{Paru dans \textit{Arabic Sciences and Philosophy}, Cambridge University Press, vol. 26(2), p.~185-214.}}
\date{}
\author{Erwan Penchèvre}

\begin{document}

\newcommand{\inputfig}[1]{\input{#1}}

\maketitle

Nous avons tenté de restituer ici les mathématiques qui
président aux théories planétaires exposées par l'astronome syrien Ibn
al-\v{S}\=a\d{t}ir (1304-1375) dans son ouvrage \textit{Nih\=ayat
  al-S\=ul}. Dans la lignée des astronomes de l'école de Mar\=agha, en
composant des mouvements de rotation à vitesse angulaire constante,
Ibn al-\v{S}\=a\d{t}ir atteint deux objectifs. Non seulement il élimine
tout recours aux excentriques et aux points équants ; mais il décrit aussi
longitudes et latitudes planétaires par une méthode unique, sans adjoindre
aucun orbe en sus des orbes nécessaires à la description
des seuls mouvements en longitude. Une meilleure compréhension des rotations
comme transformations spatiales lui permet cette grande économie
de moyens. Dans notre commentaire, nous prenons pour exemple la planète Vénus
dont les latitudes posent un problème intéressant. C'est aussi
l'occasion d'offrir l'édition critique d'un chapitre de la
\textit{Nih\=ayat al-S\=ul} consacré aux mouvements en latitude de
Mercure et Vénus.

\paragraph{Préliminaires} 
\footnote{Nous remercions Roshdi Rashed sans qui ce travail n'aurait pu
  voir le jour.}Nous n'analyserons pas en détail ce qu'est un orbe pour Ibn
al-\v{S}\=a\d{t}ir. Qu'il suffise de le concevoir comme un corps à
symétrie sphérique. Chaque orbe est doté d'un centre. Le centre de
certains orbes est aussi centre du Monde. Enfin, et surtout, Ibn
al-\v{S}\=a\d{t}ir utilise chaque orbe comme un référentiel solide en
mouvement. Pour qu'un corps solide à symétrie sphérique puisse faire
office de référentiel spatial, il faut distinguer un plan attaché à
l'orbe et passant par son centre, et une direction dans ce plan. Dans
tout ce qui suit, on représentera plus commodément chaque orbe par une
sphère, un cercle et un point~: à savoir, une sphère dont le centre
est le centre de l'orbe, un grand cercle section de cette sphère par
son plan, et un point de ce cercle situé dans la direction distinguée
par rapport au centre.

Chaque orbe est mobile par rapport à l'orbe qui le porte (c'est-à-dire
au sein du référentiel constitué par l'orbe qui le porte). Un seul
mouvement est admis pour chaque orbe~: un mouvement de rotation
uniforme autour d'un axe passant par son centre. Autrement dit, le
centre de chaque orbe est immobile au sein de l'orbe qui le porte. Il
n'y a aucune restriction quant à l'inclinaison de l'axe par rapport au
plan de l'orbe~: tel orbe a son axe perpendiculaire à son plan, tel
autre a son axe perpendiculaire au plan de l'orbe qui le porte. Mais
tous ces mouvements sont relatifs, et chaque orbe hérite aussi <<~par
accident~>> (au sens aristotélicien du terme) du mouvement de l'orbe
qui le porte, qui hérite lui-même du mouvement de l'orbe qui le porte,
et ainsi de suite. Le mouvement de chaque orbe est donc composé d'un
mouvement propre de rotation uniforme par rapport à l'orbe qui le
porte, et du mouvement absolu de l'orbe qui le porte. Seul le neuvième
orbe, porté par aucun autre, n'est animé que de son mouvement propre
de rotation uniforme.

Ce cadre conceptuel donne lieu à trois combinaisons possibles, toutes
présentes dans le texte d'Ibn al-\v{S}\=a\d{t}ir. Nous les décrirons
au moyen des trois figures (i), (ii) et (iii)
page~\pageref{fig2}. Dans chacune de ces trois figures il y a deux
orbes representés chacun par une sphère, un grand cercle de cette
sphère et un point sur ce cercle. Le premier orbe dont le point est
$P$ porte le second dont le point est $Q$. Dans les figures (i) et
(ii), le point $P$ de l'orbe portant est choisi par commodité au
centre de l'orbe porté. Au contraire, dans la figure (iii), les deux
orbes ayant même centre $O$, il a fallu choisir un autre point~: le
point $P$ est alors un point à l'intersection des plans des deux
orbes. Dans tout ce paragraphe, on désignera chaque orbe par son
point. L'orbe de $P$ porte donc l'orbe de $Q$. Ainsi, le lieu occupé
par l'orbe de $Q$ est constant au sein du référentiel solide constitué
par l'orbe de $P$. Mais l'orbe de $Q$ est animé d'un mouvement de
rotation uniforme sur lui-même, autour d'un axe représenté par un
vecteur sur chaque figure. Peu importe ici le sens de cette rotation~;
seule importe la direction de son axe, perpendiculaire au plan de
l'orbe de $P$ dans la figure (i), mais perpendiculaire au plan de
l'orbe de $Q$ dans les figures (ii) et (iii). La figure (i) présente
cette particularité que la trajectoire du point $Q$ au sein du
référentiel constitué par l'orbe de $P$ n'est pas dans le plan de
l'orbe de $Q$~: le point $Q$ décrit un petit cercle parallèle au plan
de l'orbe de $P$. Enfin, dans chacune des trois figures, l'orbe de $P$
est lui-même animé d'un mouvement de rotation uniforme autour d'un axe
représenté par un vecteur d'origine $O$. Peu importe que ce vecteur
soit perpendiculaire au plan de l'orbe de $P$ -- comme il l'est dans
ces figures -- ou non. Ce dernier mouvement entraîne le point $P$ et
aussi l'orbe de $Q$.

Chacune des planètes supérieures (Mars, Jupiter, Saturne) est portée
par un système d'orbes: un <<~orbe parécliptique~>>, un <<~orbe
incliné~>>, un <<~orbe déférent~>>, un <<~orbe rotateur~>> et un
<<~orbe de l'épicycle~>>. La figure (i) représente alors la relation
entre l'orbe incliné (orbe de $P$) et l'orbe déférent (orbe de $Q$),
ainsi que la relation entre l'orbe déférent (orbe de $P$) et l'orbe
rotateur (orbe de $Q$). Pour Vénus, les orbes portent les mêmes noms,
mais la figure (ii) représente alors la relation entre l'orbe déférent
(orbe de $P$) et l'orbe rotateur (orbe de $Q$), ainsi que la relation
entre l'orbe rotateur (orbe de $P$) et l'orbe de l'épicycle (orbe de
$Q$). Mercure compte davantage d'orbes sans toutefois introduire de
combinaison nouvelle.

Pour chaque planète, la figure (iii) représente la relation entre
l'orbe parécliptique (orbe de $P$) et l'orbe incliné (orbe de $Q$)~:
le point $P$ est alors l'un des n{\oe}uds. Mais cette figure
représente aussi la relation entre le neuvième orbe (orbe de $P$) et
le huitième, dit de l'écliptique (orbe de $Q$)~: le point $P$ est
alors le point vernal.

\begin{figure}
\begin{center}
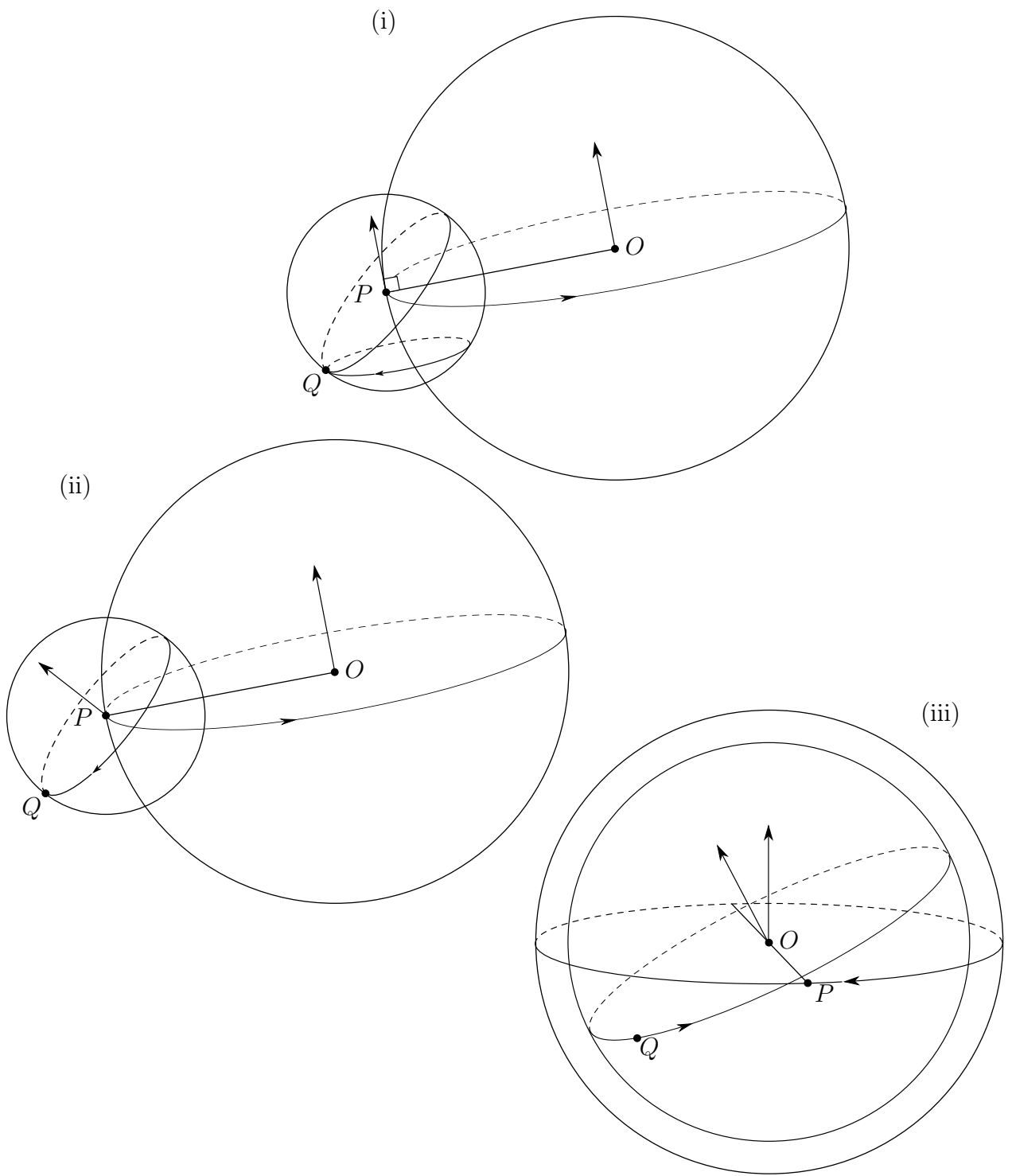
\caption{\label{fig2}Un orbe est un référentiel solide en mouvement}
\end{center}
\end{figure}

\paragraph{Précession des équinoxes}
Observons à nouveau la figure (iii). Chez Ptolémée, les étoiles fixes
sont attachées au dernier orbe et le mouvement relatif de cet orbe par
rapport aux orbes inférieurs rend compte du mouvement de précession
des équinoxes. Pour les astronomes arabes, depuis le \emph{Traité sur
  l'année solaire} certainement dû aux Banu Musa (IXème
siècle)\footnote{\textit{cf.} \cite{qurra1987}
  p.~\textit{xlvi}-\textit{lxxv} et 27-67.}, il existe un neuvième
orbe au delà des étoiles fixes situées sur le huitième orbe, et c'est
le mouvement relatif du huitième orbe par rapport au neuvième qui rend
alors compte de la précession. Tous les orbes inférieurs des planètes,
du Soleil et de la Lune héritent de ce mouvement comme les étoiles
fixes (tandis que chez Ptolémée, le Soleil en était
exempt\footnote{\textit{cf.} \cite{pedersen1974} p.~147.}).

Chez Ibn al-\v{S}\=a\d{t}ir, l'orbe parécliptique de chaque planète
est un orbe dont le centre coïncide avec le centre du Monde, comme
l'écliptique. L'orbe parécliptique reproduit, en plus petit, l'orbe de
l'écliptique et y est attaché en son centre. Le parécliptique est
animé d'un petit mouvement de rotation sur lui-même, par rapport au
référentiel constitué par l'orbe de l'écliptique~; mais comme l'axe de
rotation du parécliptique coïncide avec l'axe de l'écliptique, et que
la composée de deux rotations de même axe en est encore une, le
parécliptique est aussi en rotation uniforme autour de son axe au sein
du référentiel constitué par le neuvième orbe. Le mouvement de l'orbe
de l'écliptique par rapport au neuvième est d'un degré toutes les
soixante-dix années persanes. Le mouvement de l'orbe parécliptique par
rapport au neuvième orbe est un peu plus rapide, comme l'avait vérifié
al-Zarqalluh (XIème siècle)~: il est d'un degré toutes les soixante
années persanes, c'est-à-dire $0°1'$ par année
persane\footnote{Al-Zarqalluh (Azarquiel) avait observé que l'apogée
  solaire se meut d'un degré tous les 279 ans par rapport aux étoiles
  fixes, elles-mêmes étant animées du mouvement de précession
  (\textit{cf.} \cite{rashed1997} p.~285-287). Or
  $\dfrac{1}{70}+\dfrac{1}{279}\simeq\dfrac{1}{60}$.}.

\begin{figure}
\begin{center}
\tiny
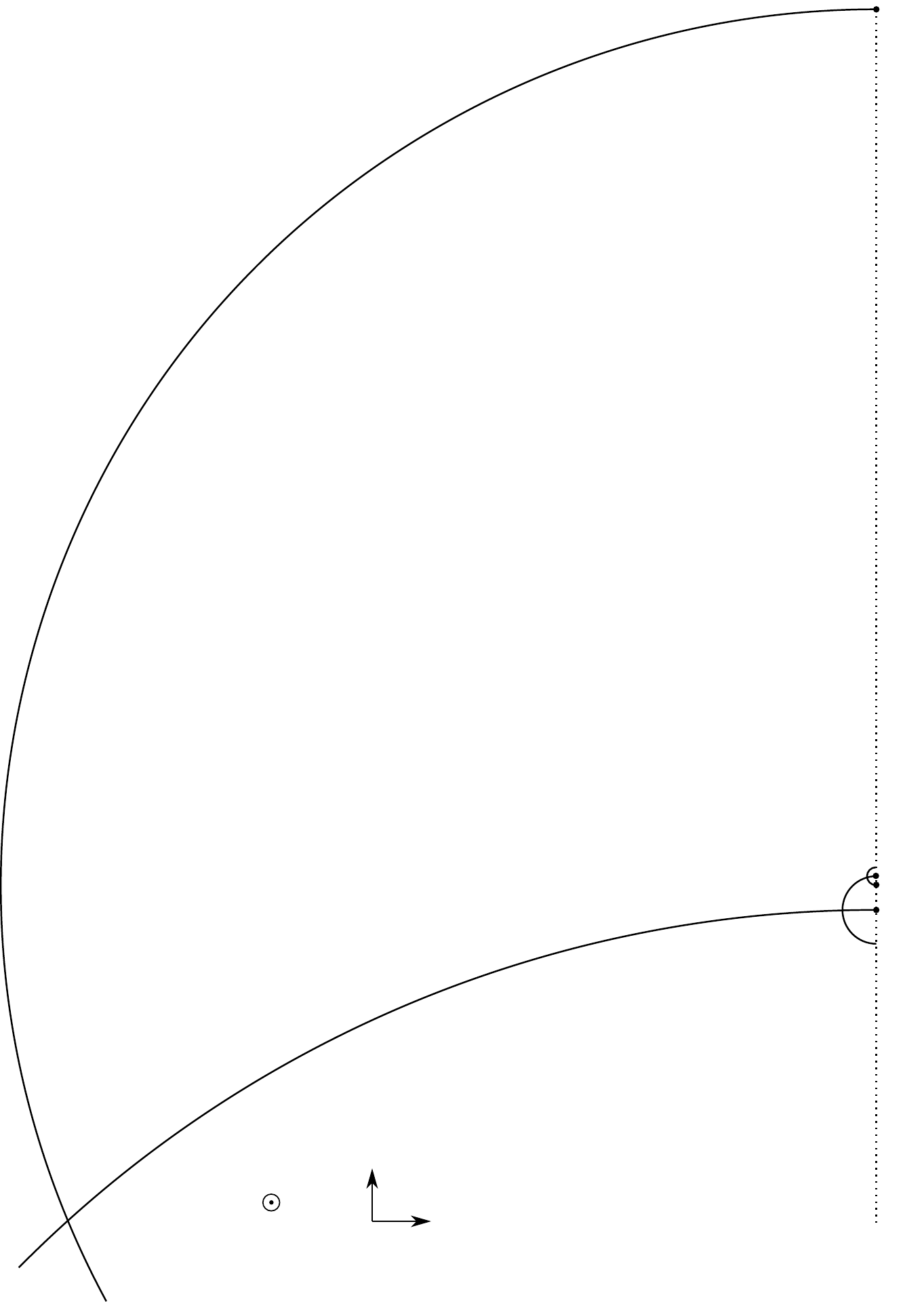
\caption{\label{fig1}Figure initiale à l'échelle}
\end{center}
\end{figure}

\begin{figure}
\begin{center}
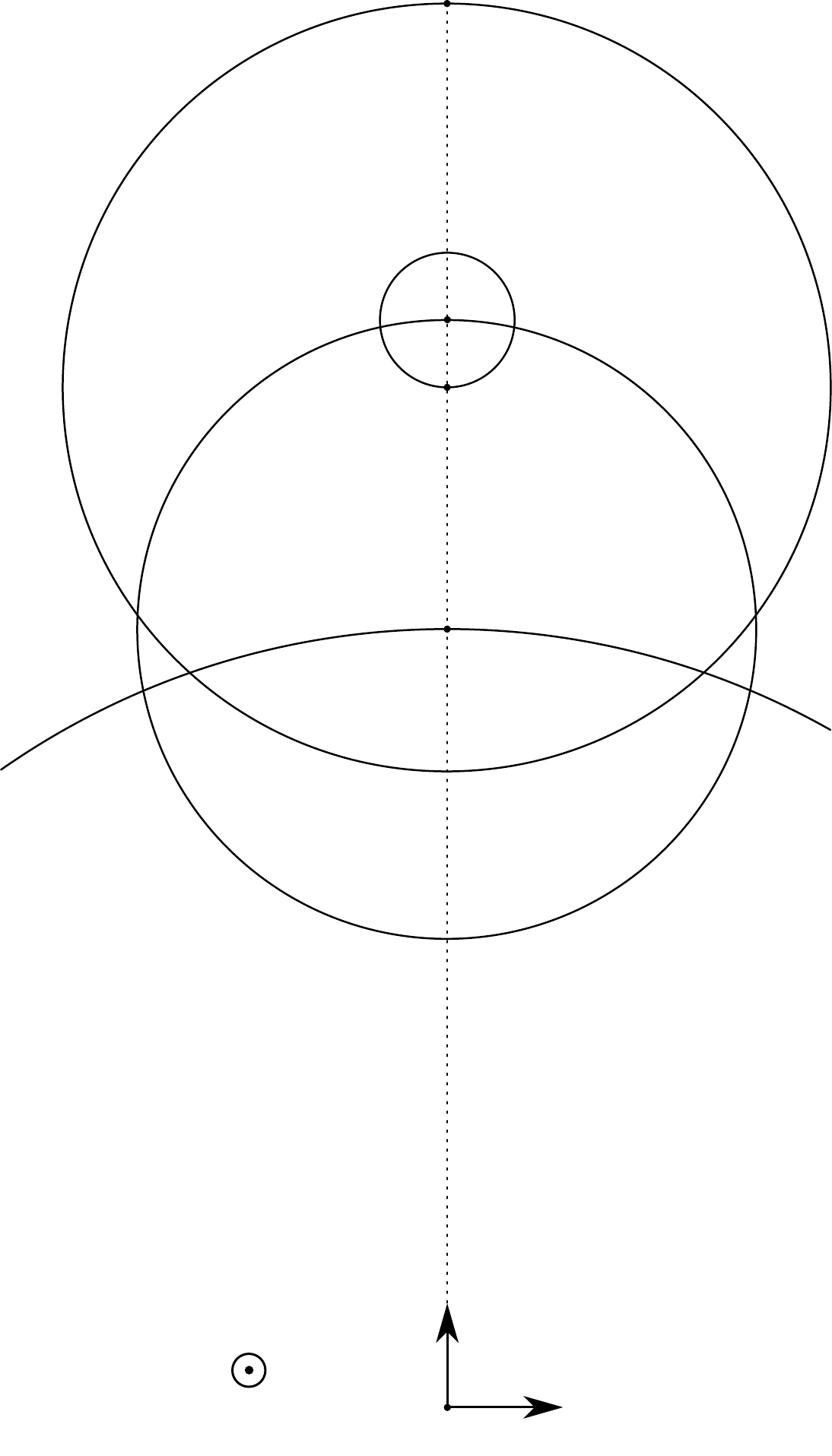
\caption{\label{fig3}Figure initiale}
\end{center}
\end{figure}

\paragraph{Figure initiale}
Soit $(\mathbf{i},\mathbf{j},\mathbf{k})$ la base canonique de
$\mathbb{R}^3$. Dans la figure initiale, les plans des orbes sont tous
rabattus dans le plan de l'écliptique $(\mathbf{i},\mathbf{j})$ par
des rotations, les centres des orbes sont tous alignés dans la
direction du vecteur $\mathbf{j}$ elle-même confondue avec la
direction du point vernal~; enfin, la ligne des n{\oe}uds, à
l'intersection des plans de l'orbe incliné et du parécliptique, est
orientée dans la direction du vecteur $\mathbf{i}$, la tête étant du
même côté que $\mathbf{i}$. Le point $O$ est le centre du Monde, $P_1$
le centre du parécliptique de Vénus, et $P_2$ le centre de son orbe
incliné. Ces trois points sont confondus $O=P_1=P_2$. Le point $P_3$
est le centre de l'orbe déférent, $P_4$ le centre de l'orbe rotateur,
$P_5$ le centre de l'orbe de l'épicycle, $P$ est le centre de Vénus,
cette planète étant elle-même un corps sphérique. La figure \ref{fig1}
précise la position de ces points~: bien qu'ils soient alignés, on
remarque que les rayons vecteurs $\overrightarrow{P_nP_{n+1}}$ ne sont
pas tous dans le même sens. Posons
$\overrightarrow{OP_3}=60\ \mathbf{j}$, alors
$$\begin{array}{l}
\overrightarrow{P_3P_4}=1;41\ \mathbf{j}\\
\overrightarrow{P_4P_5}=-0;26\ \mathbf{j}\\
\overrightarrow{P_5P}=43;33\ \mathbf{j}
\end{array}$$
où le point-virgule sépare partie entière et partie
fractionnaire\footnote{Toutes les valeurs numériques seront données en
  sexagésimal, et la virgule sert alors à séparer les rangs. Par
  exemple, $359;45,40=359+\dfrac{45}{60}+\dfrac{40}{60^2}$.}. La
figure \ref{fig1} est à l'échelle, mais il a fallu, pour ce faire,
omettre le point $O$. Par souci de lisibilité, on préfèrera utiliser
la figure \ref{fig3}, analogue mais pas à l'échelle.

Remarquons d'emblée que la figure initiale ne représente la position
des orbes à aucun instant de l'histoire passée ou future du Monde~: il
faut appliquer au point $P$ une suite de transformations géométriques
pour obtenir la position de Vénus prédite ou observée à un instant
donné.

\paragraph{Transformations géométriques}
Dans les modèles d'Ibn al-\v{S}\=a\d{t}ir, on chercherait en vain une
description mathématique du mouvement en termes de transformations ou
de différences entre deux instants. L'objet de la description
mathématique est ici l'opération transformant une figure initiale (qui
n'a d'existence qu'imaginaire) en une autre figure représentant la
configuration des astres à un instant donné. Les transformations
appliquées à Vénus $P$ sont des rotations paramétrées par trois angles
$\theta_a$, $\theta_c$, $\theta_p$. Voici la liste de ces
transformations~:
$$R(P_5,\mathbf{k},\theta_p-\theta_c),\quad R(P_5,\mathbf{j},0°30'),
\quad R(P_5,\mathbf{i},0°5'),$$
$$R(P_4,\mathbf{k},2\theta_c),\quad R(P_4,\mathbf{j},3°),
\quad R(P_4,\mathbf{i},-0°5'),$$
$$R(P_3,\mathbf{k},-\theta_c),\quad R(P_2,\mathbf{k},\theta_c),\quad
R(P_2,\mathbf{i},0°10'),\quad R(P_1,\mathbf{k},\theta_a).$$
Une rotation dont l'axe est dans la direction du vecteur $i$ ou $j$ a pour
effet d'incliner le plan d'un orbe par rapport au plan de l'orbe qui
le porte. La description d'Ibn al-\v{S}\=a\d{t}ir ne laisse aucun
doute quant à l'ordre dans lequel appliquer toutes ces
transformations.

Pour obtenir la configuration des orbes à un instant donné, il faut
donc appliquer toutes ces rotations aux points $P_3$, $P_4$, $P_5$,
$P$. L'image du point $P_3$ entraîné par les mouvements de l'orbe
parécliptique et de l'orbe incliné est~:
$$R(P_1,\mathbf{k},\theta_a)\circ R(P_2,\mathbf{i},0°10')\circ
R(P_2,\mathbf{k},\theta_c)\ (P_3)$$ Quant au point $P_4$, il est aussi
entraîné par le mouvement de l'orbe déférent et devient~:
$$R(P_1,\mathbf{k},\theta_a)\circ R(P_2,\mathbf{i},0°10')\circ
R(P_2,\mathbf{k},\theta_c)\circ R(P_3,\mathbf{k},-\theta_c)\ (P_4)$$
Le point $P_5$ est aussi entraîné par le mouvement de l'orbe rotateur
et devient~:\label{modele1}
$$\begin{array}{l} 
  R(P_1,\mathbf{k},\theta_a)\circ 
    R(P_2,\mathbf{i},0°10')\circ 
    R(P_2,\mathbf{k},\theta_c)\\
  \quad\circ
    R(P_3,\mathbf{k},-\theta_c)\circ R(P_4,\mathbf{i},-0°5')\circ
    R(P_4,\mathbf{j},3°)\circ R(P_4,\mathbf{k},2\theta_c)\ (P_5)
\end{array}$$
Enfin le point $P$, aussi entraîné par l'orbe de l'épicycle, devient~:
\label{composees}
$$\begin{array}{l}
  R(P_1,\mathbf{k},\theta_a)\circ R(P_2,\mathbf{i},0°10')\circ 
    R(P_2,\mathbf{k},\theta_c)\\
  \quad\circ R(P_3,\mathbf{k},-\theta_c)\circ R(P_4,\mathbf{i},-0°5')\circ 
    R(P_4,\mathbf{j},3°)\circ R(P_4,\mathbf{k},2\theta_c)\\
  \quad\circ R(P_5,\mathbf{i},0°5')\circ R(P_5,\mathbf{j},0°30')\circ 
    R(P_5,\mathbf{k},\theta_p-\theta_c)\ (P)
\end{array}$$

\paragraph{Trajectoire}
Dans le cadre théorique que l'on vient de décrire, l'on ne peut guère
parler de trajectoire dans l'espace. Il y a bien toutefois une
trajectoire dans l'ensemble des valeurs des paramètres~:
$$\theta_a=\omega_at+\phi_a$$
$$\theta_c=(\omega_m^{\astrosun}-\omega_a)t+(\phi_m^{\astrosun}-\phi_a)$$
$$\theta_p=\omega_pt+\phi_p$$ où l'on a distingué par un exposant les
grandeurs $(\omega_m^{\astrosun},\phi_m^{\astrosun})$ en rapport avec
le modèle du Soleil\footnote{Ainsi, le mouvement de Vénus est couplé
  au mouvement du Soleil puisque la vitesse angulaire de son orbe
  incliné par rapport au parécliptique est
  $\omega_m^{\astrosun}-\omega_a$. Le fait que le mouvement des
  planètes soit couplé au mouvement du Soleil est un défaut essentiel
  de tous les modèles pré-coperniciens.}. On a~:
$$\phi_m^{\astrosun}=280;9,0,\quad\omega_m^{\astrosun}=359;45,40
\text{ par année persane}$$
$$\phi_a=77;52,10,\quad \omega_a=0;1\text{ par année persane}$$
$$\phi_p=320;50,19,\quad \omega_p=225;1,48,41\text{ par année persane}$$

\paragraph{Unité de temps}
Le temps absolu est mesuré par le mouvement du neuvième orbe, seul
orbe dont le mouvement propre est absolu (les mouvements propres des
autres orbes sont relatifs à l'orbe qui les porte).

%La mesure terrestre d'une durée ne peut être qu'approximative. 
%Par exemple, Ibn al-\v{S}\=a\d{t}ir calcule que le neuvième orbe 
%se meut d'une minute d'arc en le temps <<~qui suffit à un homme 
%pour compter rapidement jusqu'à six~>>.
En un lieu donné, la durée entre deux culminations du Soleil (de midi
à midi du lendemain) est approximativement constante.

Considérons deux unités possibles pour la mesure du temps~:
\begin{itemize}
\item le \emph{jour sidéral}~: durée d'une rotation complète du ciel
  visible autour du pôle nord (temps de retour d'une étoile fixe en
  son lieu par rapport au pôle nord et à l'horizon)
\item le \emph{jour solaire}~: durée entre deux culminations du Soleil
  (de midi à midi du lendemain)
\end{itemize}
Le neuvième orbe étant supposé mû d'un mouvement uniforme, le jour
sidéral est de durée constante (sauf à tenir compte du mouvement de
précession, le mouvement propre du huitième orbe qui est infime). En
revanche, la durée du jour solaire varie au fil des saisons, ce n'est
donc pas une bonne unité. En effet, soit midi au Soleil. Un jour
sidéral plus tard, le Soleil aura presque effectué une rotation
complète d'est en ouest à cause du mouvement du neuvième orbe~; pas
tout à fait cependant, à cause du mouvement du parécliptique du Soleil
qui l'aura légèrement entraîné vers l'est, d'un peu moins qu'un
degré. Il ne sera donc pas encore midi au Soleil~: le jour sidéral est
plus court que le jour solaire. Cet effet se fera sentir d'autant plus
que l'angle entre le parécliptique et la trajectoire diurne du Soleil
sera faible (ainsi la différence sera plus importante lors des
solstices que lors des équinoxes). De toute façon, le mouvement du
Soleil le long de son parécliptique n'est pas uniforme.

Tout calendrier civil adoptant pourtant le jour solaire comme unité,
il est commode de choisir comme unité astronomique un <<~jour solaire
moyen~>>~: durée entre deux culminations, non plus du Soleil, mais
d'un point imaginaire mû d'un mouvement uniforme le long de l'équateur
et faisant le tour du ciel exactement en même temps que le
Soleil\footnote{c'est-à-dire en une année tropique}. C'est l'unité
adoptée par Ibn al-\v{S}\=a\d{t}ir~: \textit{al-yawm
  bi-laylatihi}. Quand une date est donnée en temps civil
(c'est-à-dire un certain nombre d'heures à compter de midi au Soleil),
il faut la corriger en lui ajoutant une certaine <<~équation du
temps~>>, \textit{ta`d{\=\i}l al-'ay\=am}, qui s'annule lors de
l'équinoxe de printemps quand le Soleil et le point imaginaire
coïncident sur l'équateur\footnote{En fait, comme on va le voir, le
  point imaginaire suit le Soleil \emph{moyen} et il est donc situé un
  peu à l'Ouest du Soleil au moment de l'équinoxe.}. L'heure est un
vingt-quatrième de jour solaire moyen.

Soit $t_0$ la date d'un équinoxe de printemps, et $t$ un instant donné
quelconque, dates exprimées en jours solaires moyens. Les mêmes dates
exprimées en temps civil (jours solaires, et heures comptées à partir
de midi au soleil) seront notées $\tau_0$ et $\tau$, et en jours
sidéraux $T_0$ et $T$.

Le point imaginaire se déplace le long de l'équateur à la vitesse du
Soleil moyen~; son ascension droite varie donc comme la longitude
$\lambda_m$ du Soleil moyen. \`A l'instant $t$, elle vaut
$\lambda_m(t)-\lambda_m(t_0)$.

De combien le Soleil s'est-il déplacé vers l'Est depuis le dernier
équinoxe de printemps~? C'est précisément son ascension droite
$\alpha(t)$.

Le jour sidéral est plus court que le jour solaire, d'autant qu'il
faut au mouvement diurne pour récupérer le déplacement du Soleil vers
l'Est. La vitesse du mouvement diurne est d'environ 15° par heure
(elle est de 360° par jour sidéral, mais le jour sidéral compte un peu
moins que 24 heures). On a donc~:
$$(T-T_0)-(\tau-\tau_0)\simeq\frac{\alpha(t)}{15°}\text{ heures}$$ 
Il ne s'agit que d'une approximation, puisque pendant cette durée, le
Soleil se déplace légèrement davantage vers l'Est...

De même, le jour sidéral est plus court que le jour solaire moyen,
d'autant qu'il faut au mouvement diurne pour récupérer le déplacement
du point imaginaire vers l'Est le long de l'équateur~:
$$(T-T_0)-(t-t_0)\simeq\frac{\lambda_m(t)-\lambda_m(t_0)}{15°}\text{ heures}$$

Ibn al-\v{S}\=a\d{t}ir donne $\lambda_m(t_0)=-2°1'7''$. L'équation du
temps est donc\footnote{On adopte ici la convention française quant au signe de l'équation du temps $E$.}~:
$$E=(t-t_0)-(\tau-\tau_0)
\simeq\frac{-\lambda_m(t)-2°1'7''+\alpha(t)}{15°}\text{ heures}$$

Chez Ibn al-\v{S}\=a\d{t}ir, l'unité de temps est le jour solaire
moyen, ou l'année persane de 365 jours solaires moyens. L'<<~époque~>>
(c'est-à-dire l'instant $t=0$) est à midi du 24 décembre 1331. Il
s'agit probablement d'une date en temps solaire moyen (les valeurs des
paramètres à cette date ayant été obtenus par le calcul à partir d'une
observation antérieure\footnote{Une remarque de Kennedy et Roberts
  (\cite{roberts1959} p.~232) au sujet du \textit{z{\=\i}j} d'Ibn
  al-\v{S}\=a\d{t}ir le confirme.}, un jour d'équinoxe de
printemps). Le temps GMT (Greenwich Mean Time) est un temps solaire
moyen mesurant comme précédemment le déplacement uniforme d'un point
imaginaire le long de l'équateur, mais ce point ne coïncide pas avec
le Soleil lors de l'équinoxe de printemps~: l'équation du temps n'a
donc pas la même origine que chez Ibn al-\v{S}\=a\d{t}ir, elle est
d'environ 8 min à l'équinoxe de printemps et elle s'annule vers le 15
avril. Il faut aussi tenir compte du décalage horaire à la longitude
de Damas, 36°18'23''. L'époque d'Ibn al-\v{S}\=a\d{t}ir est donc, en
temps GMT, le 24 décembre 1331 à 12 h 8 min
$-\dfrac{36°18'23''}{15°}\simeq$ 9 h 43 min.

Si l'on utilise la méthode d'Ibn al-\v{S}\=a\d{t}ir décrite dans cet
article pour calculer la longitude de Vénus à cette date, en
négligeant l'effet de l'opérateur $M$ (\textit{cf.} paragraphe
suivant) et en appliquant la méthode d'interpolation de Ptolémée
(\textit{cf.} p.~\pageref{ptolemee}), on trouve $264°23'$. L'Institut
de Mécanique Céleste et de Calcul des \'Ephémérides de l'Observatoire
de Paris\footnote{%
\url{http://www.imcce.fr/fr/ephemerides/formulaire/form_ephepos.php},
  choisir <<~coordonnées moyennes de la date~>>).} donne $264°21'$.

\begin{figure}
\begin{center}
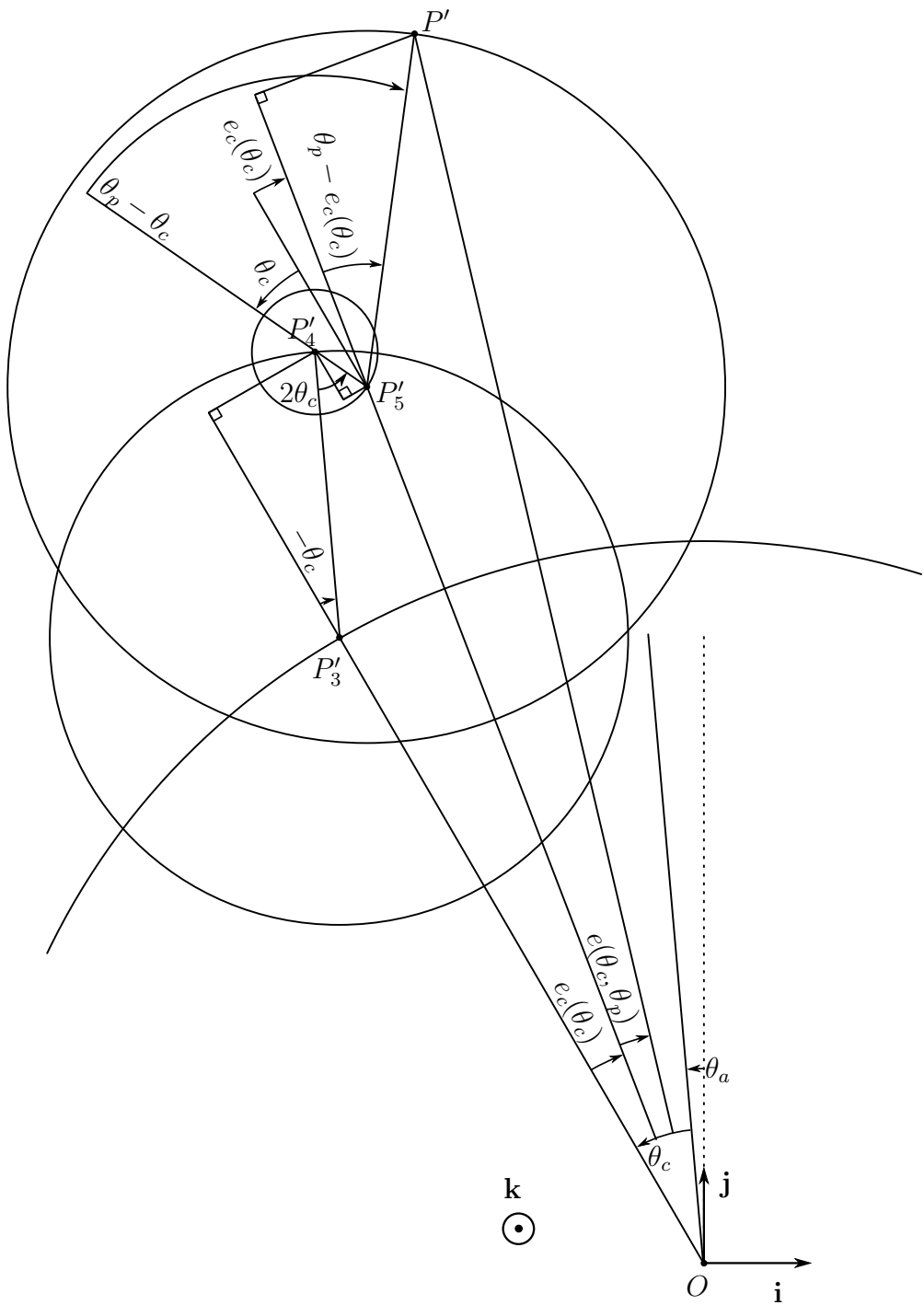
\caption{\label{fig4}Transformations planes}
\end{center}
\end{figure}

\paragraph{Transformations planes}
Si $R$ et $S$ sont deux rotations dans l'espace, il existe une
rotation $T$ telle que $R\circ S=T\circ R$, à savoir, la rotation de
même angle que $S$ et dont l'axe est l'image par $R$ de l'axe de
$S$. Appliquons cette relation de commutation à toutes les composées
de rotations décrites p.~\pageref{composees} ci-dessus, de sorte à
réécrire à droite toutes les rotations dont l'axe est dans la
direction du vecteur $\mathbf{k}$. Par exemple
$$R(P_1,\mathbf{k},\theta_a)\circ R(P_2,\mathbf{i},0°10')=
R(P_2,\mathbf{u},0°10')\circ R(P_1,\mathbf{k},\theta_a)$$ 
où le vecteur $\mathbf{u}$ est l'image de $\mathbf{i}$ par
$R(P_1,\mathbf{k},\theta_a)$. On fait de même avec les inclinaisons
des plans des petits orbes. Il existe donc une composée de rotations
$M$ telle que l'image du point $P$ soit
$$M\circ R(P_2,\mathbf{u},0°10')\ (P')$$
où
$$P'=R(P_1,\mathbf{k},\theta_a)\circ R(P_2,\mathbf{k},\theta_c)
\circ R(P_3,\mathbf{k},-\theta_c)\circ R(P_4,\mathbf{k},2\theta_c)
\circ R(P_5,\mathbf{k},\theta_p-\theta_c)\ (P)$$
On introduit de même les points $P_3'$, $P_4'$, $P_5'$ suivants 
(voir figure \ref{fig4})~:
$$P_3'=R(P_1,\mathbf{k},\theta_a)\circ R(P_2,\mathbf{k},\theta_c)\ (P_3)$$
$$P_4'=R(P_1,\mathbf{k},\theta_a)\circ R(P_2,\mathbf{k},\theta_c)\circ
R(P_3,\mathbf{k},-\theta_c)\ (P_4)$$
$$P_5'=R(P_1,\mathbf{k},\theta_a)\circ R(P_2,\mathbf{k},\theta_c)\circ
R(P_3,\mathbf{k},-\theta_c)\circ R(P_4,\mathbf{k},2\theta_c)\ (P_5)$$
Tous ces points sont dans le plan de la figure initiale~: calculer
leurs positions relève entièrement de la géométrie plane.

\paragraph{Les équations de Vénus} 
Le calcul de la position du point $P'$ revient à résoudre trois
triangles rectangles (\textit{cf}. fig.~\ref{fig4} p.~\pageref{fig4},
et chapitre 23 de la \emph{Nih\=ayat al-S\=ul}). On calcule d'abord
l'angle $e_c(\theta_c)=(OP_3',OP_5')$ appelé <<~équation du
centre~>>~:
$$e_c(\theta_c)=
-\arcsin\left(\frac{P_3P_4\sin\theta_c+P_4P_5\sin\theta_c}{OP_5'}\right)$$
où
$$OP_5'=\sqrt{(P_3P_4\sin\theta_c+P_4P_5\sin\theta_c)^2+
(OP_3+P_3P_4\cos\theta_c-P_4P_5\cos\theta_c)^2}$$
On calcule ensuite une autre <<~équation~>> $e(\theta_c,\theta_p)=
(OP_5',OP')$~:
$$e(\theta_c,\theta_p)=
\arcsin\left(\frac{P_5P\sin(\theta_p-e_c(\theta_c))}{OP'}\right)$$
où
$$OP'=\sqrt{(P_5P\sin(\theta_p-e_c(\theta_c)))^2+
(OP_5'+P_5P\cos(\theta_p-e_c(\theta_c)))^2}$$
La position de $P'$ est alors donnée par $OP'$ et par l'angle suivant~:
$$\left(\mathbf{j},\overrightarrow{OP'}\right)=
\theta_a+\theta_c+e_c(\theta_c)+e(\theta_c,\theta_p)$$
On remarque que~:
$$e_c(360°-\theta_c)=-e_c(\theta_c)$$
$$e(360°-\theta_c,\ 360°-\theta_p)=-e(\theta_c,\theta_p)$$

\begin{figure}
\begin{center}
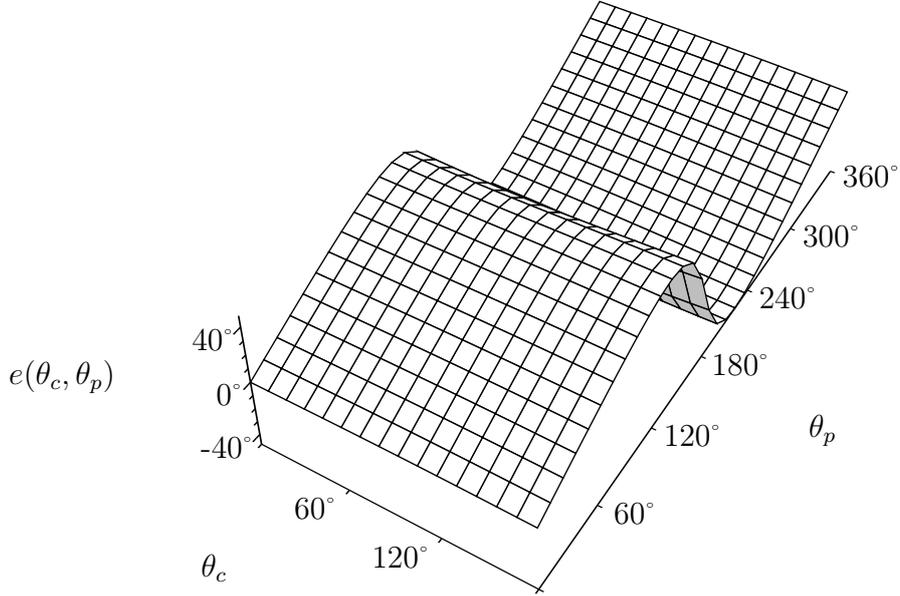
\caption{\label{fig5}Surface représentative de la fonction
  $(\theta_c,\theta_p)\mapsto e(\theta_c,\theta_p)$}
\end{center}
\end{figure}

\begin{figure}
\begin{center}
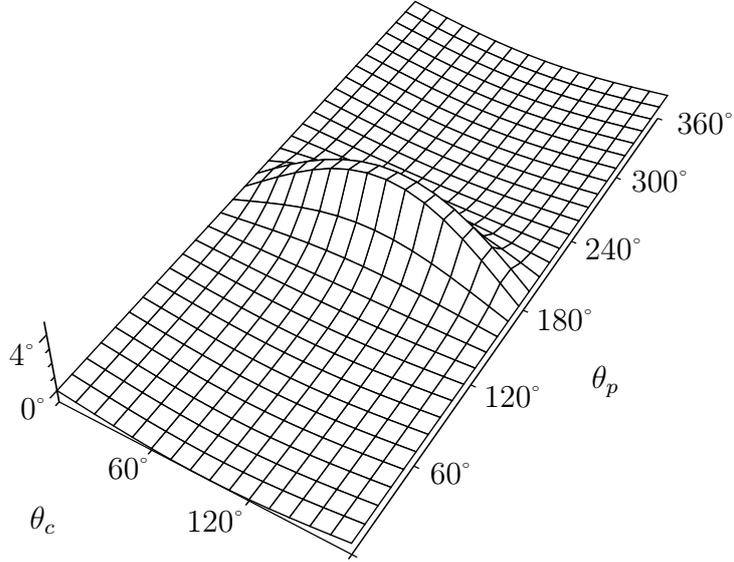
\caption{\label{fig6}$(\theta_c,\theta_p)\mapsto
  \left(e(0,\theta_p)+\chi(\theta_c)(e(180°,\theta_p)-e(0,\theta_p))\right)
  -e(\theta_c,\theta_p)$}
\end{center}
\end{figure}

\paragraph{Fonction de deux variables et interpolation}\label{ptolemee}
Les valeurs de $e_c(\theta_c)$ et $e(\theta_c,\theta_p)$ devront être
reportées dans des tables~; mais $e(\theta_c,\theta_p)$ est fonction
de deux variables, et il faudrait donc construire une table pour
chaque valeur de $\theta_c$. Les astronomes, à partir de
Ptolémée\footnote{\textit{cf.} \cite{pedersen1974}, p.~84-89.}
utilisaient dans ce contexte une méthode d'interpolation visant à
calculer une valeur approchée de $e(\theta_c,\theta_p)$ au moyen d'un
produit d'une fonction de $\theta_c$ par une fonction de
$\theta_p$. \`A $\theta_p$ donné, Ibn al-\v{S}\=a\d{t}ir interpole les
valeurs de $e$ entre $e(0,\theta_p)$ et $e(180°,\theta_p)$ au moyen de
la formule suivante (qui \emph{n'est pas} linéaire en $\theta_c$)~:
$$e(\theta_c,\theta_p)\simeq e(0,\theta_p)
+\chi(\theta_c)(e(180°,\theta_p)-e(0,\theta_p))$$
Le coefficient d'interpolation $\chi$ est défini par~:
$$\chi(\theta_c) =\frac{\max\vert e(\theta_c,\cdot)\vert-\max\vert
  e(0,\cdot)\vert} {\max\vert e(180°,\cdot)\vert-\max\vert
  e(0,\cdot)\vert}$$ 
Ce coefficient est fonction d'une seule variable,
on peut donc le calculer pour une série de valeurs de $\theta_c$ et en
faire une table. Pour ce faire, on remarque que le maximum $\max\vert
e(\theta_c,\cdot)\vert$ est atteint quand la droite $(O,P')$ est
tangente au cercle de l'orbe de l'épicycle. On a alors\footnote{On le
  montre aisément par le calcul en posant
  $z^{-1}=OP_5'+P_5P\cos(\theta_p-e_c(\theta_c))$ et en remarquant que
  la quantité à maximiser $\tan^2(e(\theta_c,\theta_p))$ est
  alors un polynôme de degré 2 en $z$.}
$$\max\vert e(\theta_c,\cdot)\vert=\arcsin\frac{P_5P}{OP_5'}$$ 
(pour mieux comprendre l'effet de cette méthode d'interpolation, voir les
figures \ref{fig5} et \ref{fig6} où l'on a tracé les surfaces
représentatives de $e(\theta_c,\theta_p)$ et de l'erreur commise en
appliquant la formule approchée ci-dessus -- les axes sont gradués en
degrés).

\begin{figure}
\begin{center}
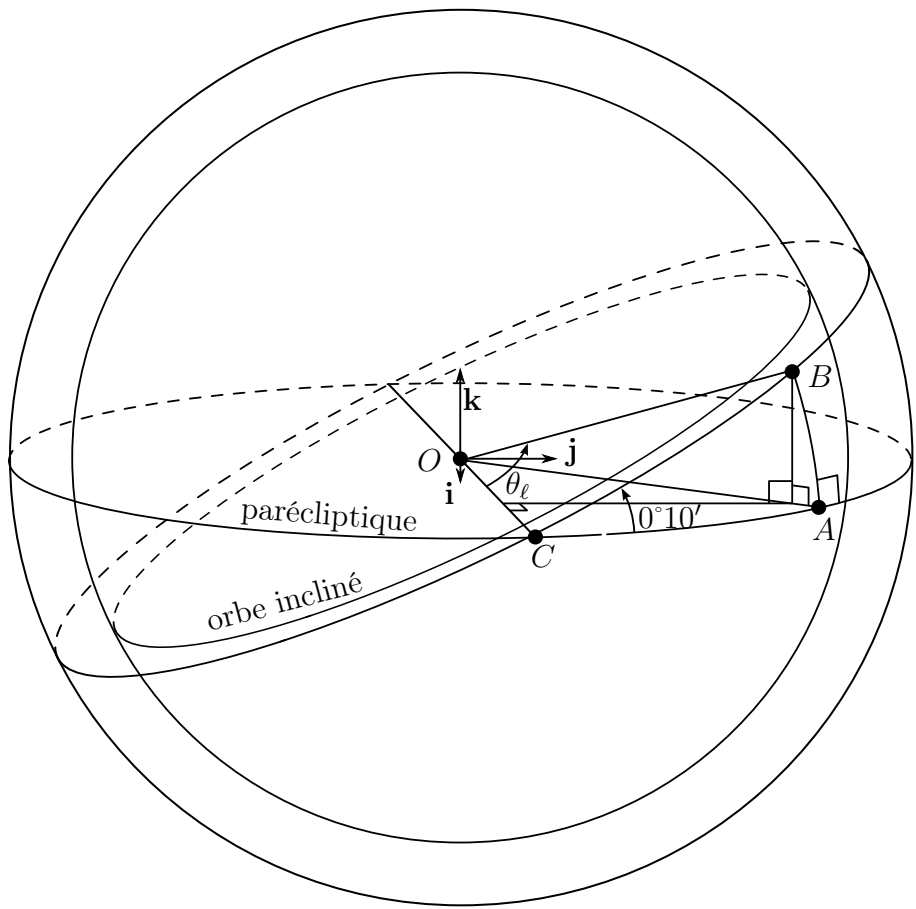
\caption{\label{fig7}Coordonnées de $R(P_2,\mathbf{u},0°10')(P')$}
\end{center}
\end{figure}

\paragraph{Trigonométrie sphérique} 
On va à présent calculer les coordonnées sphériques du point
$R(P_2,\mathbf{u},0°10')(P')$ par rapport à l'écliptique. L'angle
$\theta_\ell$ formé entre le vecteur $\mathbf{u}$ et la direction du
point $P'$ vaut
$$\theta_\ell=\theta_c+e_c(\theta_c)+e(\theta_c,\theta_p)+90°$$
L'étude du triangle sphérique $ABC$ (voir figure \ref{fig7}) donne~:
$$\tan\left(\overrightarrow{OC},\overrightarrow{OA}\right)
=\cos(0°10')\times\tan\theta_\ell$$
La longitude du point $R(P_2,\mathbf{u},0°10')(P')$ par rapport
à l'écliptique, en prenant la direction du point vernal, c'est-à-dire 
$\mathbf{j}$, comme origine, est donc (au moins quand
$\vert\theta_\ell\vert< 90°$)
$$\arctan(\cos(0°10')\times\tan\theta_\ell)-90°+\theta_a$$
et sa latitude est~:
$$\left(\overrightarrow{OA},\overrightarrow{OB}\right)=
\arcsin(\sin(0°10')\times\sin\theta_\ell)$$

\begin{figure}
\begin{center}
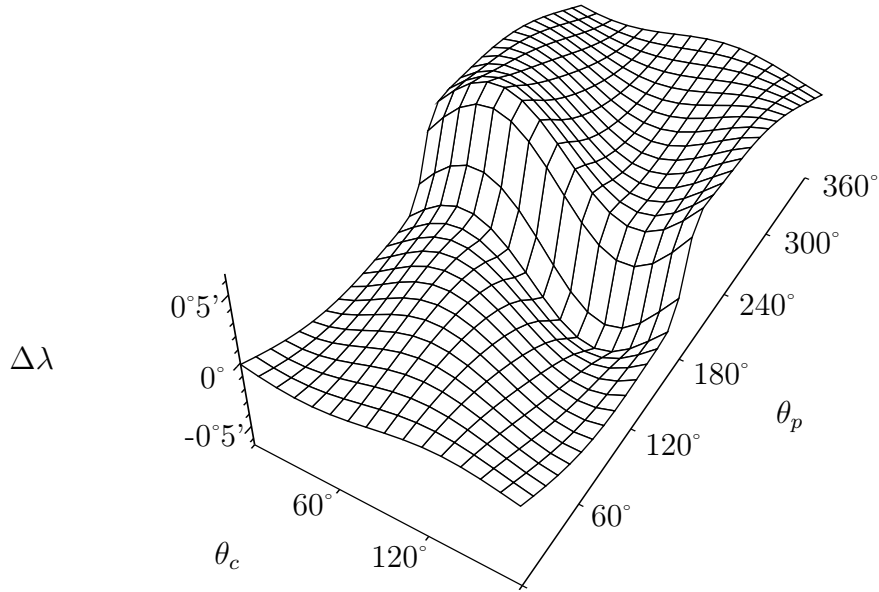
\caption{\label{fig8}Erreur en longitude si l'on néglige $M$}
\end{center}
\end{figure}

\paragraph{Les inclinaisons des petits orbes}
Les rotations qui composent $M$ ont leurs axes contenus dans le plan
de l'orbe incliné $OCB$ et sont d'angles petits ($3°$ au plus). Elles
auront peu d'effet sur la longitude de
$R(P_2,\mathbf{u},0°10')(P')$. Nous avons tracé la surface
représentative de l'erreur $\Delta\lambda(\theta_c,\theta_p)$ commise
si l'on néglige $M$ (\textit{cf}. fig.~\ref{fig8} p.~\pageref{fig8}).

\begin{figure}
\begin{center}
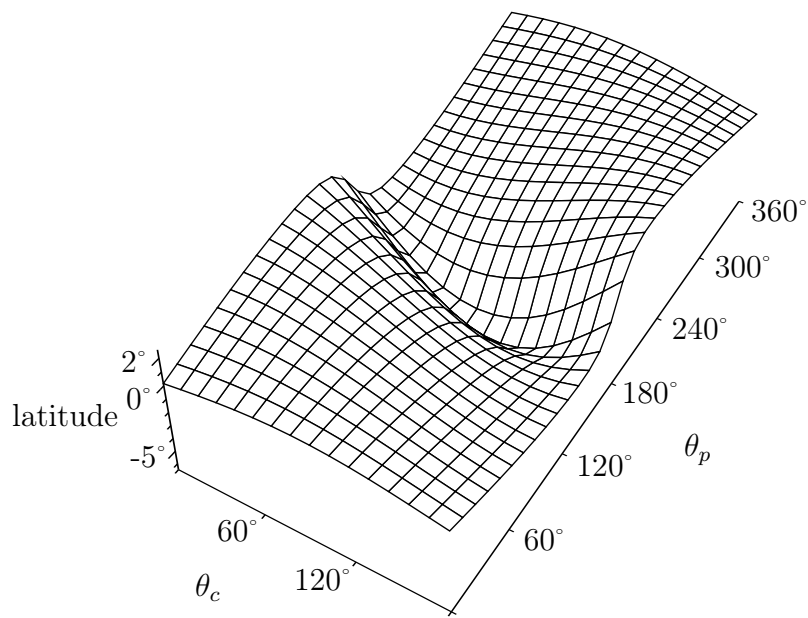
\caption{\label{fig9}Latitude de Vénus}
\end{center}
\end{figure}

En revanche, l'effet de $M$ et les inclinaisons des plans des petits
orbes ont justement été introduits dans ce modèle afin de rendre
compte des variations en latitude de Vénus. Hélas, les axes des
rotations qui composent $M$ ne passent pas par le point $O$, et de
telles rotations sont difficiles à étudier en coordonnées
sphériques\footnote{Ibn al-Hayth\=am lui-même, dans un ouvrage inégalé
  sur l'astronomie mathématique, avait renoncé à une tâche semblable~;
  \textit{cf.} \cite{rashed2006} p.~444-447.}. Ibn al-\v{S}\=a\d{t}ir
se contente ici d'une description qualitative de l'effet en latitude
de $M$. Nous donnons, fig.~\ref{fig9} p.~\pageref{fig9}, la surface
représentative de la latitude de Vénus comme fonction des deux
variables $\theta_c$ et $\theta_p$.

\begin{figure}
  \begin{center}
    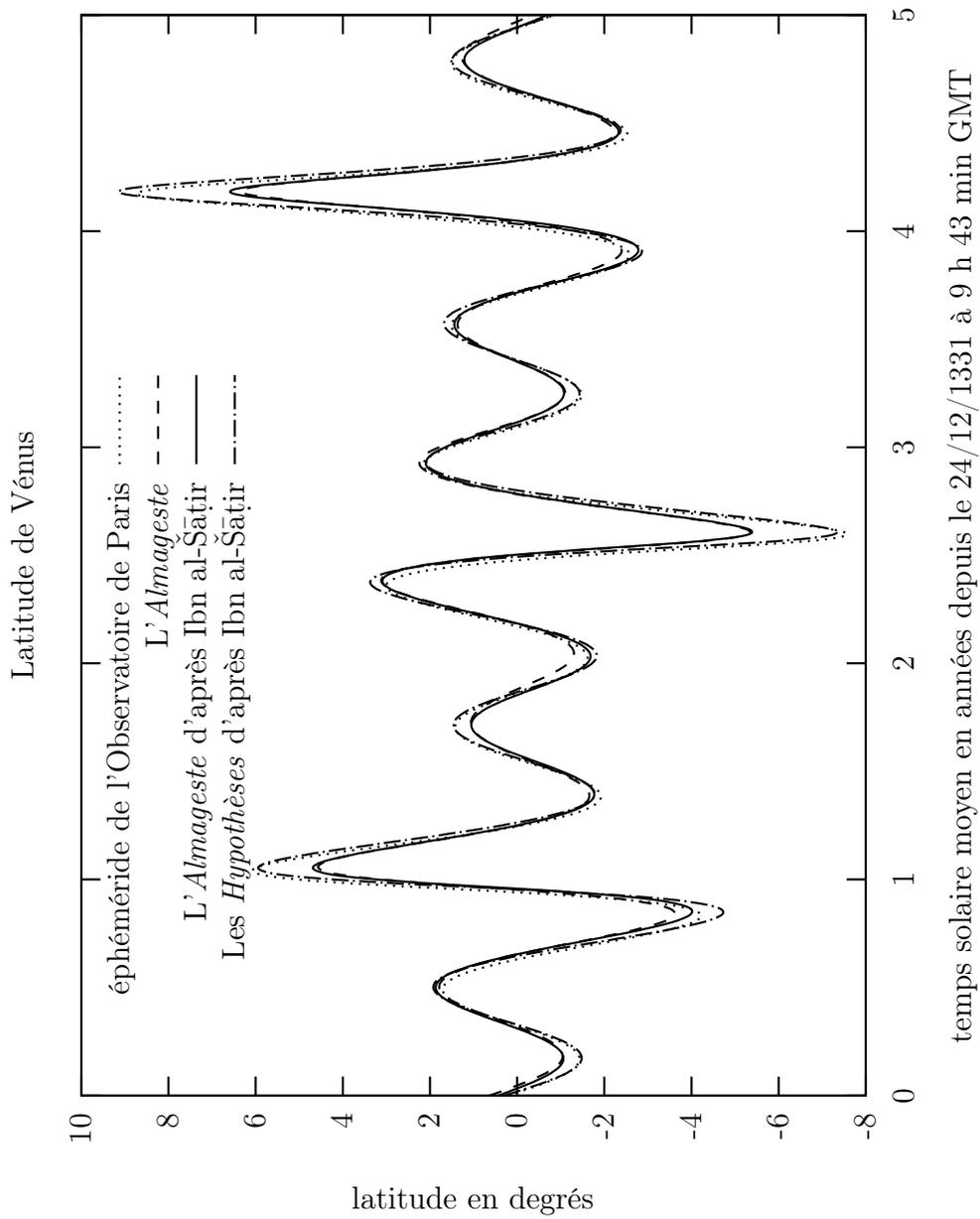
    \caption{\label{fig10}Comparaisons}
  \end{center}
\end{figure}

\paragraph{Comparaisons} Le modèle que nous avons restitué dans le
commentaire mathématique ci-dessus  p.~\pageref{modele1} est le premier
modèle décrit par Ibn al-\v{S}\=a\d{t}ir dans le chapitre vingt-cinq
de la \emph{Nih\=ayat al-S\=ul} dont nous donnons l'édition critique.
On voit que les inclinaisons des plans des orbes sont conçues dans
l'idée de reproduire les effets décrits dans l'\textit{Almageste}.
L'inclinaison de l'orbe incliné, variable dans l'\textit{Almageste}
mais constante chez Ibn al-\v{S}\=a\d{t}ir, minime, n'a que peu d'influence
sur le résultat final. Ibn al-\v{S}\=a\d{t}ir envisage ensuite une
variante plus proche de la théorie des \textit{Hypothèses planétaires}
\footnote{Sur cet autre ouvrage de Ptolémée,
  \textit{cf.} \cite{murschel1995}.}.
\`A la figure \ref{fig10}, nous donnons les latitudes de Vénus pendant cinq
ans, à partir de l'<<~époque~>> choisie par Ibn al-\v{S}\=a\d{t}ir, calculées
de quatre manières différentes : (1) par l'IMCCE de l'Observatoire de Paris,
(2) en suivant la méthode strictement ptoléméenne de l'\textit{Almageste}
mais avec les paramètres d'Ibn al-\v{S}\=a\d{t}ir à l'époque\footnote{On
  a suivi l'interprétation des méthodes de Ptolémée donnée dans \cite{pedersen1974} et \cite{swerdlow2005}.}, (3) au moyen du premier modèle d'Ibn
al-\v{S}\=a\d{t}ir censé reproduire la théorie de l'\textit{Almageste},
et (4) au moyen du second modèle d'Ibn al-\v{S}\=a\d{t}ir censé reproduire
les \textit{Hypothèses planétaires}.

\pagebreak
\section*{Chapitre vingt-cinq de la \emph{Nih\=ayat al-S\=ul}}
Des mouvements de Vénus et Mercure en latitude, en deux parties.

\paragraph{Première partie}
Des latitudes de Vénus. Chez Ptolémée, il y a deux théories à ce sujet.

\paragraph{Première théorie}
C'est celle qu'il expose dans l'\emph{Almageste} quand il dit que
l'orbe incliné de Vénus tantôt s'approche tantôt s'éloigne du plan de
l'écliptique. Aux n{\oe}uds, le centre de l'épicycle de Vénus est dans
le plan de l'écliptique, mais hors des n{\oe}uds, il est toujours au
Nord du plan de l'écliptique (au maximum, à un sixième de degré Nord,
inclinaison vers le Nord de l'orbe incliné de Vénus quand elle est à
l'Apogée\footnote{Ibn al-\v{S}\=a\d{t}ir utilise deux mots distincts
  pour désigner l'apogée~: \textit{'awj} et
  \textit{\underline{d}irwat}. Nous traduirons le premier par
  <<~Apogée~>> (avec la majuscule) et le second par <<~apogée~>>~; la
  nuance n'importe guère à la compréhension de ce
  chapitre. D'ailleurs, il n'y a qu'un seul mot pour désigner le
  périgée, \textit{\d{h}a\d{d}{\=\i}\d{d}}.}). Quand le centre de son
épicycle est aux n{\oe}uds (c'est-à-dire aux quadratures de l'Apogée),
le diamètre perpendiculaire au diamètre passant par l'Apogée de
l'épicycle et par son périgée est dans le plan de l'écliptique; quand
il est à la queue, l'Apogée de l'épicycle est incliné vers le Nord et
le périgée vers le Sud; quand il est à la tête, l'Apogée de l'épicycle
est incliné vers le Sud et son périgée vers le Nord; la grandeur de
cette inclinaison est deux degrés et demi dans l'\emph{Almageste}, et
elle s'appelle \emph{inclinaison de l'épicycle}%
\index{inclinaison!de l'épicycle}. 
Quand le centre de l'épicycle est à l'Apogée, le
diamètre passant par l'Apogée de l'épicycle et par son périgée est
dans le plan de l'orbe incliné, et le diamètre perpendiculaire est
incliné par rapport au plan de l'écliptique, de sorte que la moitié
allant de l'Apogée de l'épicycle à son périgée est au Nord de l'orbe
incliné, et l'autre moitié, au Sud; la grandeur de cette inclinaison
par rapport au plan de l'écliptique est de trois degrés et demi, et
elle s'appelle \emph{inclinaison de biais}%
\index{inclinaison!de biais}. C'est la théorie de Ptolémée dans
l'\emph{Almageste}. Certains parmi les Modernes l'ont amendée par
l'observation : ils ont trouvé que l'inclinaison de l'épicycle était
égale à l'inclinaison de biais. Les savants ont calculé les tables de
latitude de Vénus selon cette théorie.

\paragraph{Seconde théorie} 
C'est celle qu'il expose dans son livre connu sous le nom des
\emph{Hypothèses} et qui vient après le premier: l'orbe incliné de
Vénus serait d'inclinaison constante de part et d'autre. Son
inclinaison maximale est un sixième de degré vers le Nord, et de même
vers le Sud. L'inclinaison de l'épicycle est, dans les deux cas, trois
degrés et demi.

\bigskip\noindent Ni Ptolémée ni aucun autre n'avait encore pu établir
les principes du mouvement de l'épicycle et de l'orbe incliné selon la
description ci-dessus à cause de la variation d'inclinaison de
l'épicycle dans les deux cas, et de la variation des deux diamètres
dont l'un est incliné et l'autre est tantôt dans le plan de
l'écliptique, tantôt dans le plan de l'orbe incliné. On y est parvenu
grâce à Dieu.\footnote{Les caractéristiques des modèles décrits dans les
  deux paragraphes suivants sont celles que nous avons restituées dans
  notre commentaire mathématique. }

\`A cet effet, on supposera que l'inclinaison maximale de l'orbe
incliné vers le Nord est atteinte à l'Apogée de Vénus et que cette
inclinaison est un sixième de degré. C'est une inclinaison de grandeur
constante vers le Nord, et de même vers le Sud~; mais elle se meut,
elle se déplace avec l'Apogée, et les n{\oe}uds (la tête et la queue)
se déplacent aussi de par son mouvement. Supposons le centre de
l'épicycle et du rotateur sur la droite passant par
l'Apogée. Supposons le déférent dans le plan de l'orbe
incliné. Supposons l'apogée moyen du rotateur incliné vers le Sud par
rapport au plan de l'orbe incliné, d'un angle de cinq minutes d'arc,
et la moitié du rotateur commençant à l'Apogée (c'est la << première >>
moitié) inclinée vers le Nord par rapport au plan de l'orbe
incliné, formant un angle en son centre d'une grandeur de trois
degrés. Supposons l'apogée de l'épicycle incliné par rapport à
l'apogée du rotateur, vers le Nord, de cinq minutes d'arc. Le diamètre
de l'épicycle [passant par] l'Apogée reste donc dans le plan de l'orbe
incliné. Supposons la première moitié de l'épicycle (celle qui
commence à l'Apogée) inclinée d'un demi-degré vers le Nord par rapport
à la moitié du rotateur inclinée vers le Nord.

Ceci étant admis, que les orbes tournent. Que l'orbe incliné tourne
d'un quart de cercle, le déférent d'un quart de cercle aussi, et le
rotateur d'un demi-cercle dans le même temps; alors l'inclinaison de
l'épicycle s'inverse, le diamètre perpendiculaire au diamètre de
l'Apogée arrive dans le plan de l'écliptique, et l'autre diamètre est
incliné par rapport au plan de l'écliptique d'un angle d'une grandeur
de deux parts et demi; ceci, à condition que la variation
d'inclinaison de l'épicycle soit comme l'indique l'\emph{Almageste}.

Quant à l'amendement fait par les Modernes (pour qui l'épicycle est
incliné aux n{\oe}uds, à l'Apogée et au périgée de trois degrés et
demi), et la théorie qu'expose Ptolémée dans les \emph{Hypothèses},
nous supposerons l'orbe incliné comme on l'a indiqué dans la première
théorie. Le déférent est dans le plan de l'orbe incliné. \`A l'Apogée,
l'apogée du rotateur est incliné par rapport à l'orbe incliné, de cinq
minutes d'arc vers le Sud~; l'apogée de l'épicycle est incliné par
rapport à l'apogée du rotateur, de cinq minutes d'arc vers le Nord~;
la moitié du rotateur\footnote{Ici, le manuscrit porte le mot
  <<~épicycle~>>, \textit{cf.} texte arabe. C'est pourtant le plan
  du rotateur qui doit subir l'inclinaison de $3°30$ ; l'épicyle en
  hérite par composition.}
commençant à l'Apogée est inclinée vers le
Nord de trois degrés et demi. C'est la voie sur laquelle s'appuie
l'amendement fait par les Modernes.

L'inclinaison maximale de l'épicycle est, vers le Sud, $8;40$, et vers
le Nord, $1;3$, quand l'inclinaison maximale de biais\footnote{Il
  est difficile de donner un sens précis à ces deux valeurs. Peut-être
  $1°3'$ désigne-t-elle la latitude Nord maximale atteinte par
  Vénus quand le centre de l'épicycle est à la queue~: le diamètre de
  l'épicycle passant par son apogée est alors incliné de $2°30'$ par
  rapport au plan de l'écliptique, et on a
$$1°3'\simeq 2°30'\times\frac{43;33}{60+43;33}.$$
  La valeur $8°40'$ semble en revanche être erronée, bien que les
  latitudes de Vénus puissent dépasser cette valeur dans le modèle
  fondé sur l'<<~amendement fait par les Modernes~>>.} est $2;30$.
La dernière latitude
Nord est toujours, pour l'orbe incliné, un sixième de degré. La tête
de Vénus précède l'Apogée d'un quart de cercle, et sa queue le suit
d'un quart de cercle. Les dernières latitudes Nord et Sud, pour l'orbe
incliné, sont à l'Apogée et au périgée.

\paragraph{Remarque}
Concernant la latitude de Vénus, tout ce qu'ont fait les Anciens et
les Modernes en concevant des ceintures qui s'approchent ou une
inclinaison variable de l'épicycle est impossible. Qui est bien versé dans
cet art ne l'ignore pas.

\paragraph{Deuxième partie}
Des latitudes de Mercure. Chez Ptolémée, il y a deux théories, la
première dans l'\emph{Almageste} et la seconde dans les
\emph{Hypothèses}.

\paragraph{Première théorie}
Nous supposons un orbe incliné par rapport au parécliptique,
d'inclinaison maximale atteinte en l'Apogée, vers le Sud, un
demi-degré et un quart. Cette situation à l'Apogée est réciproque de
celle au périgée, où l'inclinaison est vers le Nord. De plus, la
ceinture de l'orbe incliné tend à se rapprocher de la ceinture de
l'écliptique jusqu'à ce qu'elles se confondent, puis elle s'incline à
nouveau d'autant dans l'autre sens (et de même pour le lieu situé à
l'opposé, [au périgée]). Or ceci est impossible: on ne peut le
concevoir dans l'astronomie de Ptolémée, et ni lui ni aucun autre n'a
abordé le problème du mobile de ce mouvement.

D'autre part, quand le centre de l'épicycle est au milieu de l'arc
entre les n{\oe}uds, à l'Apogée, à la dernière latitude Sud, il faut
supposer que le diamètre passant par l'Apogée et le périgée de
l'épicycle est dans le plan de l'orbe incliné, et que l'autre diamètre
est incliné de sorte que la première moitié de l'épicycle soit au Sud,
l'autre moitié au Nord, et que l'angle de cette inclinaison (dite
\emph{de biais}\index{inclinaison!de biais}) soit sept degrés. Quand
le centre de l'épicycle est aux n{\oe}uds, le diamètre perpendiculaire
à l'Apogée est dans le plan de l'écliptique, et l'autre diamètre est
incliné par rapport au plan de l'écliptique, son périgée vers le Nord
et son Apogée vers le Sud; l'angle de cette inclinaison est six degrés
et un quart, et on l'appelle \emph{inclinaison de
  l'épicycle}\index{inclinaison!de l'épicycle}. Ainsi, le centre de
l'épicycle est constamment au Sud du plan de l'écliptique, ou bien
dans ce plan (aux n{\oe}uds). C'est cette théorie qu'on trouve dans
l'\textit{Almageste}.

\paragraph{Seconde théorie}
C'est la théorie des \emph{Hypothèses}. Là, l'inclinaison de Mercure
est constante, vers le Sud à l'Apogée, d'une grandeur de six degrés
(et l'autre extrémité, au périgée, est donc inclinée vers le Nord de
six degrés). Quand le centre de l'épicycle est entre les deux
n{\oe}uds, le diamètre passant par l'Apogée coïncide avec l'orbe
incliné, et l'autre diamètre est incliné (ainsi que la deuxième moitié
de l'épicycle) vers le Nord d'un angle de six degrés et demi; on
l'appelle \emph{inclinaison de biais}. Quand le centre de l'épicycle
est [aux] n{\oe}uds, le diamètre perpendiculaire au diamètre de
l'Apogée est dans le plan de l'écliptique, et le diamètre de l'Apogée
est incliné, son périgée vers le Nord et son Apogée vers le Sud, d'un
angle égal à l'inclinaison de biais, c'est-à-dire six parts et
demi. C'est la seconde théorie; mais ni Ptolémée ni aucun autre n'a pu
établir des principes permettant de concevoir ces inclinaisons sans
perturber les mouvements en longitude. Nous avons pu concevoir ces
deux théories~; gloire à Dieu.

\bigskip\noindent Voici la théorie sur laquelle on s'est appuyé. Nous
supposons que l'inclinaison maximale de l'orbe incliné est vers le Sud
à l'Apogée et vers le Nord au périgée (cette inclinaison maximale est
de six degrés dans les \emph{Hypothèses}). Nous supposons que l'apogée
du rotateur est incliné vers le Nord de cinq minutes d'arc, que la
seconde moitié du rotateur est inclinée vers le Sud d'un angle de six
degrés, un demi-degré et un huitième de degré par rapport au plan de
l'orbe incliné, que l'apogée de l'épicycle est incliné par rapport à
l'apogée du rotateur vers le Sud de cinq minutes d'arc, et que sa
seconde moitié est inclinée de sept degrés vers le Nord par rapport au
plan de l'orbe incliné (elle est donc inclinée d'un quart et un
huitième de degré par rapport au plan du rotateur). Quand les orbes se
meuvent, l'orbe incliné d'un quart de cercle, le déférent aussi, et
l'épicycle d'un demi-cercle, alors l'inclinaison du diamètre de
l'Apogée devient six degrés et un quart, et le diamètre qui lui est
perpendiculaire passe dans le plan de l'écliptique. Or il en est ainsi
à l'observation, c'est donc la voie sur laquelle on s'est appuyé; elle
ne prend pas en considération que, selon Ptolémée, le centre de
l'épicycle serait toujours au Sud, mais il est certes revenu de cette
opinion dans les \emph{Hypothèses}. Sache cela.

La seconde théorie est conçue comme suit. Nous supposons l'orbe
incliné comme on l'a décrit ci-dessus, de sorte qu'il coupe le
parécliptique aux n{\oe}uds et que l'inclinaison maximale soit vers le
Sud à l'Apogée et vers le Nord au périgée. Nous supposons que le
centre de l'épicycle est à l'Apogée au milieu de l'arc compris entre
les n{\oe}uds. Nous supposons que l'Apogée du rotateur est incliné de
cinq minutes d'arc vers le Nord, et que sa seconde moitié est inclinée
vers le Nord, et l'autre moitié vers le Sud, d'un angle de six degrés
et demi; on appelle cela \emph{inclinaison de biais}. Nous supposons
l'apogée de l'épicycle incliné par rapport à l'apogée du rotateur de
cinq minutes d'arc vers le Sud; alors le diamètre de l'Apogée est dans
le plan de l'orbe incliné. Nous supposons que la seconde moitié de
l'épicycle est inclinée de six degrés et demi par rapport au plan de
l'orbe incliné. Que l'on se représente cela, puis que les orbes se
meuvent jusqu'à ce que le centre de l'épicycle soit au n{\oe}ud. Le
diamètre perpendiculaire au diamètre de l'Apogée coïncide avec le plan
de l'écliptique, et le diamètre de l'Apogée reste incliné de six
degrés et demi par rapport au plan de l'écliptique. La moitié [de
  l'épicycle] contenant son Apogée est vers le Sud, et l'autre moitié
vers le Nord. C'est ainsi qu'il faut concevoir les mouvements des
orbes en latitude, de façon à ne pas perturber les mouvements en
longitude.

\paragraph{Remarque}
Supposer, au besoin, que les orbes sont dans les plans des [orbes qui]
les portent ne perturbe pas les mouvements en longitude de manière
trop importante. Tu sais que ceci est analogue à l'équation du
déplacement\footnote{L'\emph{équation du déplacement} est, pour chaque
  astre, l'angle $\left(\left(\overrightarrow{OC},\overrightarrow{OA}\right)
  -\theta_\ell\right)$
  de notre figure \ref{fig7} p.~\pageref{fig7}. Le raisonnement peu
  rigoureux d'Ibn al-\v{S}\=a\d{t}ir est à peu près le suivant. Pour
  chaque astre, aussi bien que pour la Lune qu'il a traitée dans un
  chapitre antérieur, l'inclinaison de l'orbe incliné a peu
  d'influence sur le mouvement en longitude (l'équation du
  déplacement, qu'il sait calculer, est faible). Donc les inclinaisons
  des autres orbes auront \textit{a fortiori} une influence
  négligeable sur les mouvements en longitude.} de la Lune corrigée de
l'orbe incliné au parécliptique. Puisque la dernière latitude de la
Lune est de cinq degrés, et que la variation maximale due au
déplacement est de six minutes d'arc et deux tiers, l'omission de
l'équation du déplacement des trajectoires selon les latitudes à la
ceinture de l'écliptique aura pour effet une perturbation d'une minute
d'arc et un tiers, par degré de latitude, au plus: c'est une grandeur
négligeable.

Calcule le déplacement si tu le souhaites. En voici la méthode. Prends
l'élongation entre l'astre corrigé et la tête de cet astre. Prends
l'équation qui correspond à cette élongation dans la table de
l'équation du déplacement de la Lune. Multiplie cela par la latitude
de cet astre, et divise le produit par la latitude maximale de la Lune
(c'est cinq degrés). En sort l'équation qu'il faut ajouter à la
correction de l'astre si le nombre entré dans la table était compris
entre quatre-vingt-dix et cent quatre-vingt ou entre deux cent
soixante-dix et trois cent soixante degrés, et qu'il faut sinon
soustraire à la correction de l'astre; reste la vraie correction,
rapportée à l'écliptique. Cette équation est maximale dans les
octants, et nulle aux n{\oe}uds et aux dernières latitudes de part et
d'autre. Dieu est le plus grand et le plus savant.



\section*{Supplementary material}

# Édition critique

*الباب الخامس والعشرون في حركات الزهرة وعطارد في العرض وفيه فصلان

الفصل الأوّل في عروض الزهرة وفيه مذهبان لبطلميوس

الأوّل ما ذكره في المجسطي حيث قال أنّ مائل الزهرة يقرب ويبعد عن سطح البروج فلا يزال مركز تدوير الزهرة بين العقدتين في سطح فلك البروج وفي غير العقدتين شماليّ عن سطح البروج بمقدار أكثره سدس جزء ويقابل مائل الزهرة الشماليّ عند أوجها. وإذا كان مركز تدويرها في العقدتين، أعني عند تربيع الأوج، كان القطر القائم على القطر المارّ بأوج التدوير وحضيضه من سطح البروج ؛ فإذا كان في عقدة الذنب، كان ميل أوج التدوير إلى الشمال وحضيضه إلى الجنوب ؛ وإذا كان في عقدة الرأس، كان ميل أوج التدوير إلى الجنوب وحضيضه إلى الشمال، ومقدار هذا الميل في المجسطي جزأين ونصف ويعرف بميل التدوير. وإذا كان مركز التدوير في الأوج، كان القطر المارّ بأوج التدوير وحضيضه في سطح المائل والقطر القائم عليه مائل عن سطح البروج بحيث يكون النصف الّذي من الأوج إلى الحضيض من التدوير شمالياً عن المائل والآخر جنوبياً عنه، ومقدار هذا الميل عن سطح البروج ثلاثة أجزاء ونصف جزء ويعرف بميل الوراب ؛ هذا مذهبه في المجسطي. وقد صحّح بعض المتأخّرين بالرصد، فوجد ميل التدوير مثل ميل الوراب وعليه حساب جداول عرض الزهرة عند محقّقي هذا الفنّ.

ومذهبه الثاني ما ذكره في كتابه المعروف بالاقتصاص، وهو بعد الأوّل أنّ مائل الزهرة ثابت الميل في الجهتين، ونهاية ميله سدس جزءاً في الشمال ومثله في الجنوب، وميل التدوير فيه في الحالتين ثلاثة أجزاء ونصف جزء.

ولم يمكن بطلميوس ولا غيره إلى وقتنا هذا وضع أصول تحرّك التدوير

---

[1] في حركات الزهرة وعطارد في العرض وفيه فصلان : في هامش [د]، في حركات الزهرة وعطارد وفيه فصلان [ج]  [2] الفصل : ناقصة [ج]  [2] عروض : حركات [ج]  [5] ميل : مثل [د،ه،ج،ب]  [11] القطر : قطر [ب]  [12] المائل : ناقصة [ه]  [13] شمالياً : شماليّ [أ،ب،ج،ه]  [14] جنوبياً : جنوبيّ [أ،ب،ج،ه]  [15] جزء : ناقصة [أ،د،ب]  [15] وهذا : هذا [ه]  [16] بالرصد : في الرصد [ه،ب]  [18] أنّ : أمّا [ب]  [19] جزءاً : جزء [ب]  [21] أصول : أصولاً [د،ج]

والمائل بالصفة المذكورة لاختلاف ميل التدوير في الحالتين واختلاف القطرين الّذين يكون أحدهما مائلاً والآخر إمّا في سطح البروج أو في سطح المائل. وقد تيسّر بحمد الله تعالى.

وهو أن يفرض نهاية ميل المائل في الشمال عند أوج الزهرة ونهايته سدس جزء، وهو ميل ثابت المقدار في الشمال ومثله في الجنوب، وهو يتحرّك، فينتقل مع انتقال الأوج وينتقل بانتقاله العقدتين الرأس والذنب. ونفرض مركز التدوير والمدير على الخطّ المارّ بالأوج، ونفرض الحامل في سطح المائل، ونفرض ذروة المدير الوسطى مائلة إلى الجنوب عن سطح المائل زاوية مقدارها خمس دقائق والنصف من المدير الّذي ابتداؤه من الأوج وهو الأوّل شماليّ عن سطح المائل زاوية عن مركزه مقدارها ثلاثة أجزاء. ونفرض ذروة التدوير مائلة عن ذروة المدير إلى الشمال خمس دقائق، فيبقى قطر التدوير الأوج في سطح المائل. ونفرض نصف التدوير الأوّل، وهو الّذي ابتداؤه من الأوج، مائلاً عن نصف المدير المائل إلى الشمال نصف جزء إلى جهة الشمال أيضاً.

وإذا تقرّر ذلك ودارت الأفلاك، فدار المائل ربع دائرة ودار الحامل ربع دائرة أيضاً ودار المدير في ذلك الزمان نصف دائرة، فانقلب ميل التدوير، فصار القطر القائم على الأوجي في سطح البروج والقطر الآخر مائلاً عن سطح البروج زاوية مقدارها جزأين ونصف ؛ هذا على أنّ اختلاف ميل التدوير مثل ما ذكر في المجسطي.

وأمّا ما صحّحه المتأخّرون من كون ميل التدوير في العقدتين وفي الأوج والحضيض ثلاثة أجزاء ونصف جزء، وتصوّر مذهبه المذكور في الاقتصاص أن نفرض المائل على ما ذكر في المذهب الأوّل وفي سطحه الحامل في الأوج. وذروة المدير مائلة على المائل إلى الجنوب خمس دقائق وذروة التدوير مائلة عن ذروة المدير خمس دقائق إلى الشمال ونصف التدوير الّذي ابتداؤه من الأوج مائل إلى الشمال ثلاثة أجزاء ونصف.

قلت : وهذا الوجه هو المعتمد عليه وهو ما صحّحه المحققون.

---

² مائلاً : مائل [أ،ب] ⁵ يتحرّك : متحرّك [ه،ب] ⁶ بانتقاله : بانتقال [ه] ⁶ ونفرض : ونفرض ونفرض [أ]، ويفرض [د،ه،ج] ⁷ ونفرض : ويفرض [ج] ⁸ المدير : التدوير [ج،د،أ]، المدوير [ب] ¹⁰ عن : عند [ه،ب] ¹¹ الأوج : الأوجي [د]، ناقصة [ه] ¹² ونفرض : ويفرض [ه،ب] ¹² الّذي : ناقصة [ج] ¹² مائلاً : مائل [أ،ب،ج،ه] ¹⁶ الأوجي : الأوج [ج] ¹⁶ مائلاً : مائل [أ،ب،ج،ه] ¹⁷ مقدارها : مقدار [ه] ¹⁷ أنّ : ناقصة [ه] ¹⁸ ذكر : ذكره [ج] ²¹ نفرض : يفرض [د،ج] ²¹ ذكر : ذكره [ه] ²¹ سطحه : سطح [ه] ²² على : عن [د،ج] ²²⁻²³ وذروة التدوير ... خمس دقائق : في الهامش [د،ج]

27

ونهاية ميل التدوير في الجنوب ح‾م وفي الشمال ا‾ج وأكثر ميل الوراب ب‾ل. ونهاية العرض للمائل في الشمال أبداً سدس جزء. والرأس للزهرة قبل الأوج بربع دور والذنب بعد الأوج بربع دور. ونهاية العرض في الشمال والجنوب للمائل في الأوج والحضيض.

تنبيه: جميع ما تكلّفوه المتقدّمون والمتأخّرون في أمر عرض الزهرة من تقارب المنطقتين واختلاف ميل فلك التدوير محال وليس يخفى على من له اطلاع على كتب محقّقى هذا الفنّ.

الفصل الثاني في عروض عطارد، وفيها مذهبان لبطلميوس؛ الأوّل الّذي ذكره في المجسطي والثاني الّذي ذكره في الاقتصاص.

أمّا مذهبه الأوّل أن نفرض فلكاً مائلاً عن الممثّل نهاية ميله عند الأوج إلى الجنوب نصف وربع جزء، فيقاطع هذا المائل على العقدتين، ويصير نظير الجزء الّذي عند الأوج عند الحضيض مائلاً عنه إلى الشمال بما علمت. وإنّ منطقة المائل تقارب لمنطقة البروج حتّى تنطبق عليها، ثمّ يميل إلى الجهة الأخرى بذلك المقدار وكذلك الجزء المقابل له؛ وهذا محال: لا يمكن أن يتصوّر على هيئة بطلميوس ولم يتعرّض هو ولا غيره لمحرّك هذه الحركة.

ثمّ فرض إذا كان مركز التدوير في منتصف ما بين العقدتين في الأوج عند نهاية العرض الجنوبي أنّ القطر المارّ بأوج التدوير وحضيضه في سطح المائل والقطر الآخر مائل بحيث يكون النصف الأوّل من التدوير جنوبياً والنصف الآخر شمالياً وزاوية هذا الميل سبعة أجزاء، ويسمّى الوراب. وإذا كان مركز التدوير في العقدتين، كان القطر القائم على الأوج في سطح البروج والقطر الأوجي مائل عن سطح البروج، أمّا حضيضه فإلى جهة الشمال وأوجه إلى جهة الجنوب؛ وزاوية هذا الميل ستّة أجزاء وربع جزء، ويعرف بميل التدوير؛ فلا يزال مركز التدوير إمّا جنوبي عن سطح البروج وإمّا فيه عند العقدتين. هذا مذهبه في المجسطي.

ومذهبه الثاني الّذي ذكره في الاقتصاص أن مائل عطارد ثابت الميل في

---

[1] ا‾ج: ا‾ح [د،ب] [2] للمائل: المائل [ه،ب] [3] والذنب بعد الأوج بربع دور: ناقصة [ه] [3] والجنوب: في هامش [أ] [6] وليس: ليس [د،ه،ج،ب] [6] يخفى: مخفيّ [ب،د] [9] والثاني: الثاني [ج] [11] ويصير: غير مقروء [أ] [13] تقارب: يقارب [أ] [13] تنطبق: ينطبق [ه] [13] يميل: تميل [د،ج] [14] وكذلك: ولذلك [ه،ب] [15] على: على على [أ] [16] عند: ناقصة [د،ه،ج،ب] [18] جنوبياً: جنوبي [أ،ب،ج،د] [18] والنصف الآخر: والآخر [د،ج] [19] شمالياً: شماليّ [أ،ب،ج،د] [21] وأوجه إلى: وأمّا أوجه فإلى [د]، وأمّا أوجه إلى [ج]

28

الجهة الجنوبية عن الأوج ومقداره سدس جزء والطرف الآخر مائل إلى الشمال عند الحضيض سدس جزء. وإذا كان مركز التدوير فيما بين العقدتين، كان القطر المارّ بالأوج مطابقاً للمائل والقطر الآخر مائل بالنصف الثاني من التدوير إلى الشمال، ومقدار زاوية الميل ستّة أجزاء ونصف جزء، ويعرف بميل الوراب. وإذا كان مركز التدوير بين العقدتين، كان القطر القائم على الأوجي في سطح البروج والقطر الأوجي مائل، حضيضه إلى الشمال وأوجه إلى الجنوب، وزاوية هذا الميل مثل ميل الوراب؛ وذلك ستّة أجزاء ونصف. هذا مذهبه الثاني، ولم يمكن بطلميوس ولا غيره وضع أصول يمكن تصوّر هذه الميول فيها ولا تخلّ بالحركات الطولية.

وقد أمكننا نصوّر هذين المذهبين ولله الحمد.

أمّا المذهب المعتمد عليه وهو أن نفرض نهاية ميل المائل الجنوبي عند الأوج والشماليّ عند الحضيض، ونهاية هذا الميل في الاقتصاص سدس جزء؛ ونفرض ذروة المدير مائلة إلى الشمال خمس دقائق ونصف المدير الثاني مائل إلى الجنوب وزاوية الميل عن سطح المائل ستّة أجزاء ونصف وثمن وذروة التدوير مائلة عن ذروة المدير إلى الجنوب خمس دقائق ونصفه الثاني مائل إلى الشمال عن سطح المائل سبعة أجزاء، فيكون مائلاً عن سطح المدير ربع وثمن جزء. فإذا تحرّكت الأفلاك أمّا المائل فربع دور وكذلك الحامل وأمّا التدوير فنصف دائرة، فيصير ميل القطر ستّة وربع فيصير القطر الأوجي مائلاً ستّة وربع والقطر القائم عليه في سطح البروج. وهكذا كان في الرصد، وهذا الوجه هو المعتمد عليه ولا يلتفت إلى ما قاله بطلميوس بأن مركز التدوير يكون جنوبياً أبداً، فإنّه رجع عن ذلك في الاقتصاص، فاعرف ذلك.

وأمّا تصوّر المذهب الثاني أن نفرض المائل على الصفة المذكورة المتقدّمة بحيث يقاطع الممثّل على العقدتين ويكون نهاية الميل الجنوبي عند الأوج والشماليّ عند الحضيض؛ ونفرض مركز التدوير في الأوج في منتصف ما بين العقدتين؛ ونفرض أوج المدير مائلاً إلى الشمال خمس دقائق ونصف المدير

---

٣ مطابقاً : مطابق [أ،ب،د،ه] ٣ من التدوير : بالتدوير [ه،ب] ٥ الأوجي : الأوج [د،ه،ج،ب] ٧ مثل ميل : مسل [ه] ٧ أجزاء : ناقصة [أ،د،ب] ٨ الميول : غير مقروء [ه]، المدل [ب] ٩ الطولية : الطويلة [د،ب،ه] ١٠ أمكننا : أمكنا [ه،ج] ١٠ نصوّر : تصوّر [د،ج،ب] ١٠ هذين : هذا [أ،ب،د،ه] ١٢ والشماليّ : الشماليّ [ج] ١٢ ونهاية : ونهايته [ج] ١٣ الثاني : الأوّل [د] ١٤ الجنوب : الجنوبي [ب] ١٨ القطر ستّة وربع فيصير : ناقصة [ه،ب] ٢٠ يلتفت : تكترث [د]، تكترث [ج]، وكتب في هامش [ج] «يلتفت» أيضاً، غير مقروء [ب] ٢٢ نفرض : تفرض [ج] ٢٣ الممثّل : للممثّل [د،ج] ٢٤ والشماليّ : والشمال [ب] ٢٥ مائلاً : مائل [أ،ب،ه]

29

الثاني مائلاً إلى الشمال ونصفه الآخر مائلاً إلى الجنوب ومقدار زاوية هذا الميل ستّة أجزاء ونصف جزء، ويعرف بميل الوراب؛ ونفرض ذروة التدوير مائلة عن ذروة المدير إلى الجنوب خمس دقائق، فيكون القطر الأوجي في سطح المائل؛ ونفرض نصف التدوير الثاني مائلاً عن سطح المائل ستّة أجزاء ونصف.

وإذا تصوّرت ذلك، ثمّ تحرّكت الأفلاك إلى أن صار مركز التدوير في العقدة، طابق القطر القائم على الأوجي لسطح البروج، وبقي القطر الأوجي مائلاً عن سطح البروج ستّة أجزاء ونصف جزء ويكون النصف الّذي فيه أوجه جنوبياً والآخر شماليّاً. هكذا يجب تصوّر حركات الكواكب في العرض بحيث لا تخلّ بالحركات الطولية.

تنبيه مقتضى فرض الأفلاك في سطوح الحوامل في الأصل لا يخلّ بالحركات الطولية بما له كثير قدر وأنت تعلم أنّ ذلك بمثل نقل مقوم القمر من المائل إلى الممثّل. ولما كان نهاية عرض القمر خمسة أجزاء وكان نهاية الاختلاف للنقل ستّ دقائق وثلثي، كان أكثر ما يخلّ بإهمال النقل من المدائر العرضي إلى منطقة البروج أكثره لكلّ درجة من العرض دقيقة وثلث وهذا مقدار يسير.

وإن شئت فانقل، وطريق ذلك أن تأخذ بعد مقوّم الكوكب من رأس ذلك الكوكب، وهو العقدة الشمالية؛ فما حصل فخذ ما يقابله من جدول تعديل نقل القمر من التعديل، فما كان فاضربه في عرض ذلك الكوكب، فما بلغ فاقسمه على نهاية عرض القمر، وهو خمسة أجزاء؛ فما خرج فهو التعديل فزده على تقويم الكوكب، إن كان العدد الّذي دخلت به إلى الجدول أكثر من تسعين جزءاً إلى مائة وثمانين أو أكثر من مائتين وسبعين إلى ثلاثمائة وستّين؛ وإن كان غير ذلك، فنقصه من تقويم الكوكب، يبقى المقوّم المحقّق من فلك البروج وأكثره في انصاف الارباع وينعدم في العقدتين وعند نهاية العرض في الجهتين - والله تعالى أعلم.

---

¹ الثاني: الأوّل [د] ¹ مائلاً: مائل [أ،ب،ج،ه،] ² ونفرض: ويعرف [ه] ⁴ مائلاً: مائل [أ،ب،ه،] ⁶ الأوجي: القائم على الأوجي [ه] ⁶ مائلاً: مائل [أ،ب،ه،] ⁷ أوجه: الأوج أوجه [ه] ⁷ جنوبياً: جنوبي [أ،ب،ج،ه،] ⁸ شماليّاً: شماليّ [أ،ب،ج،ه،] ⁹ الطولية: الطويلة [د،ه،ب] ¹⁰ فرض: فروض [ه] ¹⁰ الأصل: الأصول [ب] ¹⁰ يخلّ: تخلّ [ج] ¹¹ كثير: كبير [د،ه،ب] ¹¹ بمثل: ممثل بمثل [ج] ¹¹ مقوم: في هامش [أ]، مقدّم [ه] ¹³ المدائر: المدار [د،ه،ج،ب] ¹³ العرضي: العرض [ب] ¹⁶ فخذ: خذ [د،ج] ¹⁷ نقل: ناقصة [ه،ب] ¹⁷ فاضربه: اضربه [ج] ¹⁸ فاقسمه: اقسمه [ج] ²⁰ إلى ثلاثمائة وستّين: ناقصة [ه] ²¹ كان: ناقصة [ه،ب] ²² وعند: أو عند [ه،ب]



% \printindex

\end{document}

%% file: art_fig_2.pdf_tex
%% Creator: Inkscape inkscape 0.48.3.1, www.inkscape.org
%% Accompanies image file 'art_fig_2.pdf' (pdf, eps, ps)
\begingroup%
  \makeatletter%
  \providecommand\color[2][]{%
    \errmessage{(Inkscape) Color is used for the text in Inkscape, but the package 'color.sty' is not loaded}%
    \renewcommand\color[2][]{}%
  }%
  \providecommand\transparent[1]{%
    \errmessage{(Inkscape) Transparency is used (non-zero) for the text in Inkscape, but the package 'transparent.sty' is not loaded}%
    \renewcommand\transparent[1]{}%
  }%
  \providecommand\rotatebox[2]{#2}%
  \ifx\svgwidth\undefined%
    \setlength{\unitlength}{467.510304bp}%
    \ifx\svgscale\undefined%
      \relax%
    \else%
      \setlength{\unitlength}{\unitlength * \real{\svgscale}}%
    \fi%
  \else%
    \setlength{\unitlength}{\svgwidth}%
  \fi%
  \global\let\svgwidth\undefined%
  \global\let\svgscale\undefined%
  \makeatother%
  \begin{picture}(1,1.17038808)%
    \put(0,0){\includegraphics[width=\unitlength]{art_fig_2.pdf}}%
    \put(0.62072055,0.92889507){\color[rgb]{0,0,0}\makebox(0,0)[lb]{\smash{$O$}}}%
    \put(0.34815203,0.88122622){\color[rgb]{0,0,0}\makebox(0,0)[lb]{\smash{$P$}}}%
    \put(0.29681626,0.79077737){\color[rgb]{0,0,0}\makebox(0,0)[lb]{\smash{$Q$}}}%
    \put(0.33969461,0.50120775){\color[rgb]{0,0,0}\makebox(0,0)[lb]{\smash{$O$}}}%
    \put(0.06712605,0.45353887){\color[rgb]{0,0,0}\makebox(0,0)[lb]{\smash{$P$}}}%
    \put(0.01579032,0.36309008){\color[rgb]{0,0,0}\makebox(0,0)[lb]{\smash{$Q$}}}%
    \put(0.77510817,0.22796177){\color[rgb]{0,0,0}\makebox(0,0)[lb]{\smash{$O$}}}%
    \put(0.81050366,0.17423217){\color[rgb]{0,0,0}\makebox(0,0)[lb]{\smash{$P$}}}%
    \put(0.63505429,0.1206466){\color[rgb]{0,0,0}\makebox(0,0)[lb]{\smash{$Q$}}}%
    \put(0.36539811,1.15778576){\color[rgb]{0,0,0}\makebox(0,0)[lb]{\smash{(i)}}}%
    \put(0.05252782,0.687616){\color[rgb]{0,0,0}\makebox(0,0)[lb]{\smash{(ii)}}}%
    \put(0.91681043,0.45771685){\color[rgb]{0,0,0}\makebox(0,0)[lb]{\smash{(iii)}}}%
  \end{picture}%
\endgroup%

%% file: art_fig_1.pdf_tex
%% Creator: Inkscape inkscape 0.48.3.1, www.inkscape.org
%% Accompanies image file 'art_fig_1.pdf' (pdf, eps, ps)
\begingroup%
  \makeatletter%
  \providecommand\color[2][]{%
    \errmessage{(Inkscape) Color is used for the text in Inkscape, but the package 'color.sty' is not loaded}%
    \renewcommand\color[2][]{}%
  }%
  \providecommand\transparent[1]{%
    \errmessage{(Inkscape) Transparency is used (non-zero) for the text in Inkscape, but the package 'transparent.sty' is not loaded}%
    \renewcommand\transparent[1]{}%
  }%
  \providecommand\rotatebox[2]{#2}%
  \ifx\svgwidth\undefined%
    \setlength{\unitlength}{392.4781333bp}%
    \ifx\svgscale\undefined%
      \relax%
    \else%
      \setlength{\unitlength}{\unitlength * \real{\svgscale}}%
    \fi%
  \else%
    \setlength{\unitlength}{\svgwidth}%
  \fi%
  \global\let\svgwidth\undefined%
  \global\let\svgscale\undefined%
  \makeatother%
  \begin{picture}(1,1.40860338)%
    \put(0,0){\includegraphics[width=\unitlength]{art_fig_1.pdf}}%
    \put(0.28504442,0.1254068){\color[rgb]{0,0,0}\makebox(0,0)[lb]{\smash{$\mathbf{k}$}}}%
    \put(0.41285907,0.12929647){\color[rgb]{0,0,0}\makebox(0,0)[lb]{\smash{$\mathbf{j}$}}}%
    \put(0.44392111,0.06363609){\color[rgb]{0,0,0}\makebox(0,0)[lb]{\smash{$\mathbf{i}$}}}%
    \put(0.96245872,1.39414261){\color[rgb]{0,0,0}\makebox(0,0)[lb]{\smash{$P$}}}%
    \put(0.31286538,0.27078013){\color[rgb]{0,0,0}\rotatebox{27.43957305}{\makebox(0,0)[lb]{\smash{orbe incliné (son centre $P_2$ est omis dans cette figure)}}}}%
    \put(0.37118748,1.21751218){\color[rgb]{0,0,0}\rotatebox{36.25205554}{\makebox(0,0)[lb]{\smash{orbe de l'épicycle dont le centre est $P_5$}}}}%
    \put(0.83910411,0.37549338){\color[rgb]{0,0,0}\makebox(0,0)[lb]{\smash{orbe déférent}}}%
    \put(0.95307063,0.41664516){\color[rgb]{0,0,0}\makebox(0,0)[lb]{\smash{$P_3$}}}%
    \put(0.95216916,0.43897808){\color[rgb]{0,0,0}\makebox(0,0)[lb]{\smash{$P_5$}}}%
    \put(0.95286822,0.45659327){\color[rgb]{0,0,0}\makebox(0,0)[lb]{\smash{$P_4$}}}%
    \put(0.8384262,0.47157208){\color[rgb]{0,0,0}\makebox(0,0)[lb]{\smash{orbe rotateur}}}%
  \end{picture}%
\endgroup%

%% file: art_fig_3.pdf_tex
%% Creator: Inkscape inkscape 0.48.3.1, www.inkscape.org
%% Accompanies image file 'art_fig_3.pdf' (pdf, eps, ps)
\begingroup%
  \makeatletter%
  \providecommand\color[2][]{%
    \errmessage{(Inkscape) Color is used for the text in Inkscape, but the package 'color.sty' is not loaded}%
    \renewcommand\color[2][]{}%
  }%
  \providecommand\transparent[1]{%
    \errmessage{(Inkscape) Transparency is used (non-zero) for the text in Inkscape, but the package 'transparent.sty' is not loaded}%
    \renewcommand\transparent[1]{}%
  }%
  \providecommand\rotatebox[2]{#2}%
  \ifx\svgwidth\undefined%
    \setlength{\unitlength}{314.9765435bp}%
    \ifx\svgscale\undefined%
      \relax%
    \else%
      \setlength{\unitlength}{\unitlength * \real{\svgscale}}%
    \fi%
  \else%
    \setlength{\unitlength}{\svgwidth}%
  \fi%
  \global\let\svgwidth\undefined%
  \global\let\svgscale\undefined%
  \makeatother%
  \begin{picture}(1,1.75082137)%
    \put(0,0){\includegraphics[width=\unitlength]{art_fig_3.pdf}}%
    \put(0.28044758,0.1423225){\color[rgb]{0,0,0}\makebox(0,0)[lb]{\smash{$\mathbf{k}$}}}%
    \put(0.55964302,0.15081907){\color[rgb]{0,0,0}\makebox(0,0)[lb]{\smash{$\mathbf{j}$}}}%
    \put(0.62749421,0.00739182){\color[rgb]{0,0,0}\makebox(0,0)[lb]{\smash{$\mathbf{i}$}}}%
    \put(0.54823182,0.95251661){\color[rgb]{0,0,0}\makebox(0,0)[lb]{\smash{$P_3$}}}%
    \put(0.55141216,1.37391223){\color[rgb]{0,0,0}\makebox(0,0)[lb]{\smash{$P_4$}}}%
    \put(0.54982195,1.24828857){\color[rgb]{0,0,0}\makebox(0,0)[lb]{\smash{$P_5$}}}%
    \put(0.51483818,0.01590498){\color[rgb]{0,0,0}\makebox(0,0)[lb]{\smash{$O$}}}%
    \put(0.54597346,1.70290432){\color[rgb]{0,0,0}\makebox(0,0)[lb]{\smash{$P$}}}%
  \end{picture}%
\endgroup%

%% file: art_fig_4.pdf_tex
%% Creator: Inkscape inkscape 0.48.3.1, www.inkscape.org
%% Accompanies image file 'art_fig_4.pdf' (pdf, eps, ps)
\begingroup%
  \makeatletter%
  \providecommand\color[2][]{%
    \errmessage{(Inkscape) Color is used for the text in Inkscape, but the package 'color.sty' is not loaded}%
    \renewcommand\color[2][]{}%
  }%
  \providecommand\transparent[1]{%
    \errmessage{(Inkscape) Transparency is used (non-zero) for the text in Inkscape, but the package 'transparent.sty' is not loaded}%
    \renewcommand\transparent[1]{}%
  }%
  \providecommand\rotatebox[2]{#2}%
  \ifx\svgwidth\undefined%
    \setlength{\unitlength}{371.04792567bp}%
    \ifx\svgscale\undefined%
      \relax%
    \else%
      \setlength{\unitlength}{\unitlength * \real{\svgscale}}%
    \fi%
  \else%
    \setlength{\unitlength}{\svgwidth}%
  \fi%
  \global\let\svgwidth\undefined%
  \global\let\svgscale\undefined%
  \makeatother%
  \begin{picture}(1,1.42089072)%
    \put(0,0){\includegraphics[width=\unitlength]{art_fig_4.pdf}}%
    \put(0.54307963,0.11924656){\color[rgb]{0,0,0}\makebox(0,0)[lb]{\smash{$\mathbf{k}$}}}%
    \put(0.7800841,0.12645876){\color[rgb]{0,0,0}\makebox(0,0)[lb]{\smash{$\mathbf{j}$}}}%
    \put(0.83768187,0.0047061){\color[rgb]{0,0,0}\makebox(0,0)[lb]{\smash{$\mathbf{i}$}}}%
    \put(0.74205006,0.01193265){\color[rgb]{0,0,0}\makebox(0,0)[lb]{\smash{$O$}}}%
    \put(0.7002372,0.15456289){\color[rgb]{0,0,0}\makebox(0,0)[lb]{\smash{$\theta_c$}}}%
    \put(0.76285782,0.25213178){\color[rgb]{0,0,0}\makebox(0,0)[lb]{\smash{$\theta_a$}}}%
    \put(0.6327512,0.40826778){\color[rgb]{0,0,0}\rotatebox{-68.41033099}{\makebox(0,0)[lb]{\smash{$e(\theta_c,\theta_p)$}}}}%
    \put(0.59653513,0.35319767){\color[rgb]{0,0,0}\rotatebox{-60.17891929}{\makebox(0,0)[lb]{\smash{$e_c(\theta_c)$}}}}%
    \put(0.30875815,0.85219477){\color[rgb]{0,0,0}\rotatebox{-60.67750088}{\makebox(0,0)[lb]{\smash{$-\theta_c$}}}}%
    \put(0.33336426,0.6915671){\color[rgb]{0,0,0}\makebox(0,0)[lb]{\smash{$P_3'$}}}%
    \put(0.09970124,1.22880497){\color[rgb]{0,0,0}\rotatebox{-34.88097497}{\makebox(0,0)[lb]{\smash{$\theta_p-\theta_c$}}}}%
    \put(0.29796726,0.99518288){\color[rgb]{0,0,0}\makebox(0,0)[lb]{\smash{$2\theta_c$}}}%
    \put(0.30568374,1.06573231){\color[rgb]{0,0,0}\makebox(0,0)[lb]{\smash{$P_4'$}}}%
    \put(0.33498945,1.29159872){\color[rgb]{0,0,0}\rotatebox{-79.2885516}{\makebox(0,0)[lb]{\smash{$\theta_p-e_c(\theta_c)$}}}}%
    \put(0.45267715,1.4041038){\color[rgb]{0,0,0}\makebox(0,0)[lb]{\smash{$P'$}}}%
    \put(0.2349425,1.30262209){\color[rgb]{0,0,0}\rotatebox{-60.17891929}{\makebox(0,0)[lb]{\smash{$e_c(\theta_c)$}}}}%
    \put(0.26858323,1.14449003){\color[rgb]{0,0,0}\rotatebox{-53.10747075}{\makebox(0,0)[lb]{\smash{$\theta_c$}}}}%
    \put(0.40274855,0.99540679){\color[rgb]{0,0,0}\makebox(0,0)[lb]{\smash{$P_5'$}}}%
  \end{picture}%
\endgroup%

%% file: art_fig_5.pdf_tex
%% Creator: Inkscape inkscape 0.48.3.1, www.inkscape.org
%% Accompanies image file 'art_fig_5.pdf' (pdf, eps, ps)
\begingroup%
  \makeatletter%
  \providecommand\color[2][]{%
    \errmessage{(Inkscape) Color is used for the text in Inkscape, but the package 'color.sty' is not loaded}%
    \renewcommand\color[2][]{}%
  }%
  \providecommand\transparent[1]{%
    \errmessage{(Inkscape) Transparency is used (non-zero) for the text in Inkscape, but the package 'transparent.sty' is not loaded}%
    \renewcommand\transparent[1]{}%
  }%
  \providecommand\rotatebox[2]{#2}%
  \ifx\svgwidth\undefined%
    \setlength{\unitlength}{343.4875026bp}%
    \ifx\svgscale\undefined%
      \relax%
    \else%
      \setlength{\unitlength}{\unitlength * \real{\svgscale}}%
    \fi%
  \else%
    \setlength{\unitlength}{\svgwidth}%
  \fi%
  \global\let\svgwidth\undefined%
  \global\let\svgscale\undefined%
  \makeatother%
  \begin{picture}(1,0.66948764)%
    \put(0,0){\includegraphics[width=\unitlength]{art_fig_5.pdf}}%
    \put(0.20837846,0.03960669){\color[rgb]{0,0,0}\makebox(0,0)[lb]{\smash{$\theta_c$}}}%
    \put(0.86830842,0.19212966){\color[rgb]{0,0,0}\makebox(0,0)[lb]{\smash{$\theta_p$}}}%
    \put(0.30964761,0.10278773){\color[rgb]{0,0,0}\makebox(0,0)[lb]{\smash{60°}}}%
    \put(0.39374813,0.04967195){\color[rgb]{0,0,0}\makebox(0,0)[lb]{\smash{120°}}}%
    \put(0.65657816,0.09762368){\color[rgb]{0,0,0}\makebox(0,0)[lb]{\smash{60°}}}%
    \put(0.71065532,0.18277057){\color[rgb]{0,0,0}\makebox(0,0)[lb]{\smash{120°}}}%
    \put(0.76277822,0.25973435){\color[rgb]{0,0,0}\makebox(0,0)[lb]{\smash{180°}}}%
    \put(0.814275,0.33286565){\color[rgb]{0,0,0}\makebox(0,0)[lb]{\smash{240°}}}%
    \put(0.86160441,0.40083668){\color[rgb]{0,0,0}\makebox(0,0)[lb]{\smash{300°}}}%
    \put(0.90491778,0.46622735){\color[rgb]{0,0,0}\makebox(0,0)[lb]{\smash{360°}}}%
    \put(-0.00030023,0.24976802){\color[rgb]{0,0,0}\makebox(0,0)[lb]{\smash{$e(\theta_c,\theta_p)$}}}%
    \put(0.20817406,0.16947991){\color[rgb]{0,0,0}\makebox(0,0)[lb]{\smash{-40°}}}%
    \put(0.2248223,0.22675409){\color[rgb]{0,0,0}\makebox(0,0)[lb]{\smash{0°}}}%
    \put(0.19842559,0.28728825){\color[rgb]{0,0,0}\makebox(0,0)[lb]{\smash{40°}}}%
  \end{picture}%
\endgroup%

%% file: art_fig_6.pdf_tex
%% Creator: Inkscape inkscape 0.48.3.1, www.inkscape.org
%% Accompanies image file 'art_fig_6.pdf' (pdf, eps, ps)
\begingroup%
  \makeatletter%
  \providecommand\color[2][]{%
    \errmessage{(Inkscape) Color is used for the text in Inkscape, but the package 'color.sty' is not loaded}%
    \renewcommand\color[2][]{}%
  }%
  \providecommand\transparent[1]{%
    \errmessage{(Inkscape) Transparency is used (non-zero) for the text in Inkscape, but the package 'transparent.sty' is not loaded}%
    \renewcommand\transparent[1]{}%
  }%
  \providecommand\rotatebox[2]{#2}%
  \ifx\svgwidth\undefined%
    \setlength{\unitlength}{275.23466484bp}%
    \ifx\svgscale\undefined%
      \relax%
    \else%
      \setlength{\unitlength}{\unitlength * \real{\svgscale}}%
    \fi%
  \else%
    \setlength{\unitlength}{\svgwidth}%
  \fi%
  \global\let\svgwidth\undefined%
  \global\let\svgscale\undefined%
  \makeatother%
  \begin{picture}(1,0.80018514)%
    \put(0,0){\includegraphics[width=\unitlength]{art_fig_6.pdf}}%
    \put(0.05308803,0.07918021){\color[rgb]{0,0,0}\makebox(0,0)[lb]{\smash{$\theta_c$}}}%
    \put(0.81442098,0.27141218){\color[rgb]{0,0,0}\makebox(0,0)[lb]{\smash{$\theta_p$}}}%
    \put(0.158721,0.14671128){\color[rgb]{0,0,0}\makebox(0,0)[lb]{\smash{60°}}}%
    \put(0.26556308,0.074765){\color[rgb]{0,0,0}\makebox(0,0)[lb]{\smash{120°}}}%
    \put(0.59545616,0.14026672){\color[rgb]{0,0,0}\makebox(0,0)[lb]{\smash{60°}}}%
    \put(0.66860222,0.24464204){\color[rgb]{0,0,0}\makebox(0,0)[lb]{\smash{120°}}}%
    \put(0.73930933,0.35012262){\color[rgb]{0,0,0}\makebox(0,0)[lb]{\smash{180°}}}%
    \put(0.80546259,0.44704792){\color[rgb]{0,0,0}\makebox(0,0)[lb]{\smash{240°}}}%
    \put(0.86830132,0.53941947){\color[rgb]{0,0,0}\makebox(0,0)[lb]{\smash{300°}}}%
    \put(0.92990072,0.61913951){\color[rgb]{0,0,0}\makebox(0,0)[lb]{\smash{360°}}}%
    \put(0.04007523,0.2290706){\color[rgb]{0,0,0}\makebox(0,0)[lb]{\smash{0°}}}%
    \put(0.02489851,0.30264881){\color[rgb]{0,0,0}\makebox(0,0)[lb]{\smash{4°}}}%
  \end{picture}%
\endgroup%

%% file: art_fig_7.pdf_tex
%% Creator: Inkscape inkscape 0.48.3.1, www.inkscape.org
%% Accompanies image file 'art_fig_7.pdf' (pdf, eps, ps)
\begingroup%
  \makeatletter%
  \providecommand\color[2][]{%
    \errmessage{(Inkscape) Color is used for the text in Inkscape, but the package 'color.sty' is not loaded}%
    \renewcommand\color[2][]{}%
  }%
  \providecommand\transparent[1]{%
    \errmessage{(Inkscape) Transparency is used (non-zero) for the text in Inkscape, but the package 'transparent.sty' is not loaded}%
    \renewcommand\transparent[1]{}%
  }%
  \providecommand\rotatebox[2]{#2}%
  \ifx\svgwidth\undefined%
    \setlength{\unitlength}{337.66860384bp}%
    \ifx\svgscale\undefined%
      \relax%
    \else%
      \setlength{\unitlength}{\unitlength * \real{\svgscale}}%
    \fi%
  \else%
    \setlength{\unitlength}{\svgwidth}%
  \fi%
  \global\let\svgwidth\undefined%
  \global\let\svgscale\undefined%
  \makeatother%
  \begin{picture}(1,0.99883682)%
    \put(0,0){\includegraphics[width=\unitlength]{art_fig_7.pdf}}%
    \put(0.57756517,0.37806942){\color[rgb]{0,0,0}\makebox(0,0)[lb]{\smash{$C$}}}%
    \put(0.451781,0.48503334){\color[rgb]{0,0,0}\makebox(0,0)[lb]{\smash{$O$}}}%
    \put(0.8880979,0.40797462){\color[rgb]{0,0,0}\makebox(0,0)[lb]{\smash{$A$}}}%
    \put(0.88310603,0.57854407){\color[rgb]{0,0,0}\makebox(0,0)[lb]{\smash{$B$}}}%
    \put(0.549295,0.45862495){\color[rgb]{0,0,0}\makebox(0,0)[lb]{\smash{$\theta_\ell$}}}%
    \put(0.48195155,0.44774323){\color[rgb]{0,0,0}\makebox(0,0)[lb]{\smash{$\mathbf{i}$}}}%
    \put(0.61733364,0.49681912){\color[rgb]{0,0,0}\makebox(0,0)[lb]{\smash{$\mathbf{j}$}}}%
    \put(0.50225889,0.5526642){\color[rgb]{0,0,0}\makebox(0,0)[lb]{\smash{$\mathbf{k}$}}}%
    \put(0.22301271,0.30603218){\color[rgb]{0,0,0}\rotatebox{14.93132055}{\makebox(0,0)[lb]{\smash{orbe incliné}}}}%
    \put(0.25665935,0.43242972){\color[rgb]{0,0,0}\rotatebox{-4.38866748}{\makebox(0,0)[lb]{\smash{parécliptique}}}}%
    \put(0.69602451,0.42151285){\color[rgb]{0,0,0}\makebox(0,0)[lb]{\smash{$0°10'$}}}%
  \end{picture}%
\endgroup%

%% file: art_fig_8.pdf_tex
%% Creator: Inkscape inkscape 0.48.3.1, www.inkscape.org
%% Accompanies image file 'art_fig_8.pdf' (pdf, eps, ps)
\begingroup%
  \makeatletter%
  \providecommand\color[2][]{%
    \errmessage{(Inkscape) Color is used for the text in Inkscape, but the package 'color.sty' is not loaded}%
    \renewcommand\color[2][]{}%
  }%
  \providecommand\transparent[1]{%
    \errmessage{(Inkscape) Transparency is used (non-zero) for the text in Inkscape, but the package 'transparent.sty' is not loaded}%
    \renewcommand\transparent[1]{}%
  }%
  \providecommand\rotatebox[2]{#2}%
  \ifx\svgwidth\undefined%
    \setlength{\unitlength}{331.09534141bp}%
    \ifx\svgscale\undefined%
      \relax%
    \else%
      \setlength{\unitlength}{\unitlength * \real{\svgscale}}%
    \fi%
  \else%
    \setlength{\unitlength}{\svgwidth}%
  \fi%
  \global\let\svgwidth\undefined%
  \global\let\svgscale\undefined%
  \makeatother%
  \begin{picture}(1,0.69929398)%
    \put(0,0){\includegraphics[width=\unitlength]{art_fig_8.pdf}}%
    \put(0.23049465,0.06435127){\color[rgb]{0,0,0}\makebox(0,0)[lb]{\smash{$\theta_c$}}}%
    \put(0.8633795,0.22415086){\color[rgb]{0,0,0}\makebox(0,0)[lb]{\smash{$\theta_p$}}}%
    \put(0.31830579,0.12048885){\color[rgb]{0,0,0}\makebox(0,0)[lb]{\smash{60°}}}%
    \put(0.40712201,0.06068097){\color[rgb]{0,0,0}\makebox(0,0)[lb]{\smash{120°}}}%
    \put(0.65470087,0.10729148){\color[rgb]{0,0,0}\makebox(0,0)[lb]{\smash{60°}}}%
    \put(0.71550604,0.18778502){\color[rgb]{0,0,0}\makebox(0,0)[lb]{\smash{120°}}}%
    \put(0.76957979,0.26919741){\color[rgb]{0,0,0}\makebox(0,0)[lb]{\smash{180°}}}%
    \put(0.82143589,0.34663394){\color[rgb]{0,0,0}\makebox(0,0)[lb]{\smash{240°}}}%
    \put(0.87210483,0.41714892){\color[rgb]{0,0,0}\makebox(0,0)[lb]{\smash{300°}}}%
    \put(0.91547141,0.48185098){\color[rgb]{0,0,0}\makebox(0,0)[lb]{\smash{360°}}}%
    \put(0.19725953,0.19787594){\color[rgb]{0,0,0}\makebox(0,0)[lb]{\smash{-0°5'}}}%
    \put(0.21450778,0.26451694){\color[rgb]{0,0,0}\makebox(0,0)[lb]{\smash{0°}}}%
    \put(0.18079528,0.34448623){\color[rgb]{0,0,0}\makebox(0,0)[lb]{\smash{0°5'}}}%
    \put(-0.00031147,0.28568496){\color[rgb]{0,0,0}\makebox(0,0)[lb]{\smash{$\Delta\lambda$}}}%
  \end{picture}%
\endgroup%

%% file: art_fig_9.pdf_tex
%% Creator: Inkscape inkscape 0.48.3.1, www.inkscape.org
%% Accompanies image file 'art_fig_9.pdf' (pdf, eps, ps)
\begingroup%
  \makeatletter%
  \providecommand\color[2][]{%
    \errmessage{(Inkscape) Color is used for the text in Inkscape, but the package 'color.sty' is not loaded}%
    \renewcommand\color[2][]{}%
  }%
  \providecommand\transparent[1]{%
    \errmessage{(Inkscape) Transparency is used (non-zero) for the text in Inkscape, but the package 'transparent.sty' is not loaded}%
    \renewcommand\transparent[1]{}%
  }%
  \providecommand\rotatebox[2]{#2}%
  \ifx\svgwidth\undefined%
    \setlength{\unitlength}{295.30944098bp}%
    \ifx\svgscale\undefined%
      \relax%
    \else%
      \setlength{\unitlength}{\unitlength * \real{\svgscale}}%
    \fi%
  \else%
    \setlength{\unitlength}{\svgwidth}%
  \fi%
  \global\let\svgwidth\undefined%
  \global\let\svgscale\undefined%
  \makeatother%
  \begin{picture}(1,0.80555767)%
    \put(0,0){\includegraphics[width=\unitlength]{art_fig_9.pdf}}%
    \put(0.22923626,0.06191525){\color[rgb]{0,0,0}\makebox(0,0)[lb]{\smash{$\theta_c$}}}%
    \put(0.84682369,0.2398365){\color[rgb]{0,0,0}\makebox(0,0)[lb]{\smash{$\theta_p$}}}%
    \put(0.25683047,0.13231442){\color[rgb]{0,0,0}\makebox(0,0)[lb]{\smash{60°}}}%
    \put(0.35392328,0.06898829){\color[rgb]{0,0,0}\makebox(0,0)[lb]{\smash{120°}}}%
    \put(0.6439352,0.12000403){\color[rgb]{0,0,0}\makebox(0,0)[lb]{\smash{60°}}}%
    \put(0.71086572,0.2152244){\color[rgb]{0,0,0}\makebox(0,0)[lb]{\smash{120°}}}%
    \put(0.76776278,0.30774544){\color[rgb]{0,0,0}\makebox(0,0)[lb]{\smash{180°}}}%
    \put(0.83087534,0.38835025){\color[rgb]{0,0,0}\makebox(0,0)[lb]{\smash{240°}}}%
    \put(0.88519822,0.46741031){\color[rgb]{0,0,0}\makebox(0,0)[lb]{\smash{300°}}}%
    \put(0.9313337,0.53870991){\color[rgb]{0,0,0}\makebox(0,0)[lb]{\smash{360°}}}%
    \put(0.14722128,0.22653809){\color[rgb]{0,0,0}\makebox(0,0)[lb]{\smash{-5°}}}%
    \put(0.14200803,0.30032232){\color[rgb]{0,0,0}\makebox(0,0)[lb]{\smash{0°}}}%
    \put(0.13559905,0.34709476){\color[rgb]{0,0,0}\makebox(0,0)[lb]{\smash{2°}}}%
    \put(-0.0021905,0.27774926){\color[rgb]{0,0,0}\makebox(0,0)[lb]{\smash{latitude}}}%
  \end{picture}%
\endgroup%

%% file: art_fig_10.pdf_tex
%% Creator: Inkscape inkscape 0.48.5, www.inkscape.org
%% Accompanies image file 'art_fig_10.pdf' (pdf, eps, ps)
\begingroup%
  \makeatletter%
  \providecommand\color[2][]{%
    \errmessage{(Inkscape) Color is used for the text in Inkscape, but the package 'color.sty' is not loaded}%
    \renewcommand\color[2][]{}%
  }%
  \providecommand\transparent[1]{%
    \errmessage{(Inkscape) Transparency is used (non-zero) for the text in Inkscape, but the package 'transparent.sty' is not loaded}%
    \renewcommand\transparent[1]{}%
  }%
  \providecommand\rotatebox[2]{#2}%
  \ifx\svgwidth\undefined%
    \setlength{\unitlength}{384bp}%
    \ifx\svgscale\undefined%
      \relax%
    \else%
      \setlength{\unitlength}{\unitlength * \real{\svgscale}}%
    \fi%
  \else%
    \setlength{\unitlength}{\svgwidth}%
  \fi%
  \global\let\svgwidth\undefined%
  \global\let\svgscale\undefined%
  \makeatother%
  \begin{picture}(1,1.25)%
    \put(0,0){\includegraphics[width=\unitlength]{art_fig_10.pdf}}%
    \put(0.23458333,0.73770827){\rotatebox{90}{\makebox(0,0)[rb]{\smash{L'\textit{Almageste} d'après Ibn al-\v{S}\=a\d{t}ir}}}}%
    \put(0.889375,0.11520827){\rotatebox{90}{\makebox(0,0)[rb]{\smash{-8}}}}%
    \put(0.80416666,0.11520827){\rotatebox{90}{\makebox(0,0)[rb]{\smash{-6}}}}%
    \put(0.71895833,0.11520827){\rotatebox{90}{\makebox(0,0)[rb]{\smash{-4}}}}%
    \put(0.63354166,0.11520827){\rotatebox{90}{\makebox(0,0)[rb]{\smash{-2}}}}%
    \put(0.54833333,0.11520827){\rotatebox{90}{\makebox(0,0)[rb]{\smash{0}}}}%
    \put(0.463125,0.11520827){\rotatebox{90}{\makebox(0,0)[rb]{\smash{2}}}}%
    \put(0.37791666,0.11520827){\rotatebox{90}{\makebox(0,0)[rb]{\smash{4}}}}%
    \put(0.2925,0.11520827){\rotatebox{90}{\makebox(0,0)[rb]{\smash{6}}}}%
    \put(0.20729166,0.11520827){\rotatebox{90}{\makebox(0,0)[rb]{\smash{8}}}}%
    \put(0.12208333,0.11520827){\rotatebox{90}{\makebox(0,0)[rb]{\smash{10}}}}%
    \put(0.926875,0.13249994){\rotatebox{90}{\makebox(0,0)[b]{\smash{0}}}}%
    \put(0.926875,0.34562494){\rotatebox{90}{\makebox(0,0)[b]{\smash{1}}}}%
    \put(0.926875,0.55874994){\rotatebox{90}{\makebox(0,0)[b]{\smash{2}}}}%
    \put(0.926875,0.7716666){\rotatebox{90}{\makebox(0,0)[b]{\smash{3}}}}%
    \put(0.15958333,0.73770827){\rotatebox{90}{\makebox(0,0)[rb]{\smash{éphéméride de l'Observatoire de Paris}}}}%
    \put(0.926875,0.9847916){\rotatebox{90}{\makebox(0,0)[b]{\smash{4}}}}%
    \put(0.926875,1.1979166){\rotatebox{90}{\makebox(0,0)[b]{\smash{5}}}}%
    \put(0.19708333,0.73770827){\rotatebox{90}{\makebox(0,0)[rb]{\smash{L'\textit{Almageste}}}}}%
    \put(0.49566487,0.02075169){\makebox(0,0)[b]{\smash{latitude en degrés}}}%
    \put(0.983125,0.66520827){\rotatebox{90}{\makebox(0,0)[b]{\smash{temps solaire moyen en années depuis le 24/12/1331 à 9 h 43 min GMT}}}}%
    \put(0.06583333,0.66520827){\rotatebox{90}{\makebox(0,0)[b]{\smash{Latitude de Vénus}}}}%
    \put(0.27208333,0.73770827){\rotatebox{90}{\makebox(0,0)[rb]{\smash{Les \textit{Hypothèses} d'après Ibn al-\v{S}\=a\d{t}ir}}}}%
  \end{picture}%
\endgroup%

%% file: article.bbl
\begin{thebibliography}{}

\bibitem{altusi1993} Nas{\=\i}r al-D{\=\i}n
  al-T\=us{\=\i}. \textit{Nas{\=\i}r al-D{\=\i}n al-T\=us{\=\i}'s
    Memoir on astronomy: al-tadhkira f{\=\i} `ilm al-hay'a}, édition,
  traduction et commentaire par F. J. Ragep. Springer, New York, 1993.

\bibitem{alurdi1990} Mu'ayyid al-D{\=\i}n
  al-`Urd{\=\i}. \textit{Kit\=ab al-Hay'a}, édition et présentation
  par G. Saliba. Beyrouth 1990.

\bibitem{ghanem1976} E. S. Kennedy et Imad Ghanem. \textit{Ibn
  al-Sh\=a\d{t}ir, An Arab Astronomer of the Fourteenth
  Century}. Aleppo Univ. Publications, 1976.

\bibitem{hartner1974} W. Hartner. Ptolemy, Azarquiel, Ibn
  al-Sh\=a\d{t}ir and Copernicus on Mercury. A Study of
  Parameters. \textit{Archives Internationales d'Histoire des
    Sciences}, 24: 5-25, 1974.

\bibitem{hasnawi2001} Ahmad Hasnawi. La définition du mouvement dans
  la \textit{Physique} du \textit{\v{S}if\=a'}
  d'Avicenne. \textit{Arabic Sciences and Philosophy}, 11: 219-255,
  2001.

\bibitem{kennedy1966} E. S. Kennedy. Late Medieval Planetary
  Theory. \textit{Isis}, 57: 365-378, 1966 (\cite{ghanem1976}, p.~93).

\bibitem{murschel1995} Andrea Murschel. The structure and function of Ptolemy's
  physical hypotheses of planetary motion. \textit{Journal for the History
    of Astronomy}, 26: 33-61, 1995.

\bibitem{pedersen1974} Olaf Pedersen. \textit{A Survey of the
  Almagest}, Odense Univ. Press, 1974.

\bibitem{qurra1987} Th\=abit Ibn Qurra. \textit{{\OE}uvres
  d'astronomie}, édition, traduction et commentaire par Régis
  Morelon. Paris, Les Belles Lettres, 1987.

\bibitem{rashed1997} Roshdi Rashed, Régis Morelon \textit{et
  alii}. \textit{Histoire des sciences arabes. 1. Astronomie,
  théorique et appliquée.} Seuil, Paris, 1997.

\bibitem{rashed2006} Roshdi Rashed. \textit{Les mathématiques
  infinitésimales du IXème au XIème siècles}. Volume V: \textit{Ibn
  al-Haytham. Astronomie, géométrie sphérique et trigonométrie.}
  Al-Furq\=an Islamic Heritage Foundation, Londres, 2006.

\bibitem{roberts1959} E. S. Kennedy et Victor Roberts. The Planetary
  Theory of Ibn al-Sh\=a\d{t}ir. \textit{Isis}, 50: 227-235, 1959
  (\cite{ghanem1976}, p.~60).

\bibitem{roberts1957} Victor Roberts. The Solar and Lunar Theory of
  Ibn ash-Sh\=a\d{t}ir, A Pre-Copernican Copernican
  Model. \textit{Isis}, 48: 428-432, 1957 (\cite{ghanem1976}, p.~44).

\bibitem{roberts1966} Victor Roberts. The Planetary Theory of Ibn
  al-Sh\=a\d{t}ir: Latitudes of the Planets. \textit{Isis}, 57:
  208-219, 1966 (\cite{ghanem1976}, p.~81).

\bibitem{swerdlow1984} N. M. Swerdlow et
  O. Neugebauer. \textit{Mathematical astronomy in Copernicus' De
    Revolutionibus}. Springer-Verlag 1984.

\bibitem{swerdlow2005} N. M. Swerdlow. Ptolemy's Theories of the Latitude
  of the Planets in the \textit{Almagest}, \textit{Handy Tables}, and
  \textit{Planetary Hypotheses}, dans \textit{Wrong for the Right Reasons},
  ed. Jed Z. Buchwald et Allan Franklin, Springer-Verlag 2005.

\bibitem{wiedemann1928} Eilhard Wiedemann. Beiträge zur Geschichte der
  Naturwissenschaften. LXXIX. Ibn al Schâ\d{t}ir, ein arabischer
  Astronom aus dem 14. Jahrhundert. \textit{Sitzungsberichten der
    Physikalisch-Medezinischen Sozietät zu Erlangen}, 60: 317-326,
  1928 (\cite{ghanem1976}, p.~17).

\end{thebibliography}
